\documentclass[12pt]{article}
\usepackage{latexsym,amssymb,
amsthm,amsfonts,amsmath, epsfig,psfrag,array }
\makeatletter
\@addtoreset{figure}{section}

\@addtoreset{table}{section}

\makeatother

\begin{document}

\newtheorem{lemma}{Lemma}[section]
\newtheorem{cor}{Corollary}[section]
\newtheorem*{theorem}{Main Theorem}
\renewcommand{\proofname}{Proof}
\newtheorem{prop}{Proposition}[section]

\def \H{{\mathbb H}}
\def \N{{\mathbb N}}
\def \R{{\mathbb R}}
\def \E{{\mathbb E}}
\def \Z{{\mathbb Z}}
\def \S{{\mathbb S}}
\def \L{{\mathcal L}}
\def \D{{G}}
\def \l{\langle }
\def \r{\rangle }
\def \lf{\left\lfloor}
\def \rf{\right\rfloor}
\def \d{D\,}
\def \Gr{\mathop{\mathrm{Gr}}}
\def \o{\overline}
\def \wt{\widetilde}
\def \wh{\widehat}
\def \a{\alpha}
\def \b{\beta}
\def \dim{\mathop{\mathrm{dim}}}
\def \conv{\mathop{\mathrm{conv}}}

\begin{center}
{\LARGE  \bf
 Compact hyperbolic Coxeter $n$-polytopes with $n+3$ facets
}

\vspace{19pt}

{\Large Pavel~Tumarkin\footnote{Partially supported by grants of President of Russia 
MK-6290.2006.1 and NS-5666.2006.1,  INTAS grant YSF-06-10000014-5916, and RFBR grant 07-01-00390-a}
}
\vspace{3pt}

pasha@mccme.ru

\vspace{10pt}

{\small Independent University of Moscow, \
B. Vlassievskii 11, 119002 Moscow, Russia.}\\

\end{center}

\bigskip

\begin{center}
\parbox{13cm}{\small {\it Abstract.\ }
We use methods of combinatorics of polytopes together 
with geometrical and computational ones to obtain the 
complete list of  compact hyperbolic
Coxeter $n$-polytopes with $n+3$ facets, $4\le n\le 7$.
Combined with results of Esselmann~\cite{Ess2} this gives the
classification of all compact hyperbolic Coxeter $n$-polytopes
with $n+3$ facets, $n\ge 4$. Polytopes in dimensions $2$ and $3$
were classified by Poincar\'e~\cite{Po} and
Andreev~\cite{A1}.}

\end{center}

\vspace{8pt}

\section{Introduction}

A polytope in the hyperbolic space $\H^n$ is called a
{\it Coxeter polytope} if its dihedral angles are all
 integer submultiples of $\pi$. Any Coxeter polytope $P$ is a
fundamental domain of the discrete group generated by
reflections in the facets of $P$.

There is no complete classification of compact hyperbolic Coxeter polytopes.
Vinberg~\cite{V2} proved there are no such polytopes in $\H^n, n\ge
30$. Examples are known only for $n\le 8$
(see~\cite{Bu1},~\cite{Bu2}).

In dimensions $2$ and $3$ compact Coxeter polytopes were completely classified by
Poinca\-r\'e~\cite{Po} and Andreev~\cite{A1}.
Compact polytopes of the simplest combinatorial
type, the simplices, were
classified by Lann\'er~\cite{L}. Kaplinskaja~\cite{Kap} (see
also~\cite{V1}) listed  simplicial prisms, Esselmann~\cite{Ess}
classified the remaining compact $n$-polytopes with $n+2$ facets.

In the paper~\cite{ImH} Im Hof classified polytopes that can be
described by Napier cycles. These polytopes have at most $n+3$ facets.
Concerning polytopes with $n+3$ facets, Esselmann proved
the following theorem (\cite[Th.~5.1]{Ess2}):

{\it
Let $P$ be a compact hyperbolic Coxeter $n$-polytope bounded by $n+3$
facets. Then $n\le 8$; if $n=8$, then $P$ is the polytope found by
Bugaenko in~\cite{Bu2}. This polytope has the following Coxeter diagram:

\begin{center}
\epsfig{file=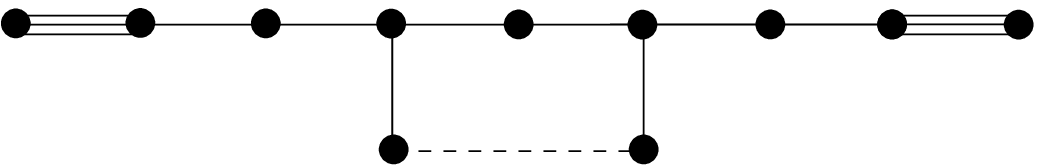,width=0.5\linewidth}
\end{center}
}

In this paper, we expand the technique derived by Esselmann
in~\cite{Ess2} and~\cite{Ess} to complete the classification of
compact hyperbolic Coxeter $n$-polytopes with $n+3$ facets.
The aim is to prove the following theorem:

\begin{theorem}
Tables~\ref{t4}--\ref{t8} contain all Coxeter diagrams of
compact hyperbolic Coxeter $n$-polytopes with $n+3$ facets for
$n\ge 4$.

\end{theorem}

The paper is organized as follows. In Section~\ref{hyp-gale} we recall basic
definitions and list some well-known properties of hyperbolic
Coxeter polytopes. We also emphasize the connection between
combinatorics (Gale diagram) and metric properties (Coxeter
diagram) of hyperbolic Coxeter polytope. In Section~\ref{tech} we
recall some technical tools from~\cite{V2} and~\cite{Ess2}
concerning Coxeter diagrams and Gale diagrams, and introduce
notation suitable for investigating of large number of diagrams.
Section~\ref{proof} is devoted to the proof of the main theorem.
The most part of the proof is computational: we restrict the
number of Coxeter diagrams in consideration, and use a computer
check after that. The bulk is to find an upper bound for the
number of diagrams, and then to reduce the number to make the
computation short enough.

\smallskip

This paper is a completely rewritten part of my Ph.D. thesis
(2004) with several errors corrected. I am grateful to my advisor
Prof.~E.~B.~Vinberg for his help. I am also grateful to
Prof.~R.~Kellerhals who brought the papers of F.~Esselmann and 
L.~Schlettwein to my attention.

\section{$\!$Hyperbolic Coxeter polytopes and Gale diagrams}
\label{hyp-gale}

In this section we list essential facts concerning hyperbolic
Coxeter polytopes, Gale diagrams of simple polytopes, and Coxeter diagrams
we use in this paper. Proofs, details and definitions in general
case may be found in~\cite{gale} and~\cite{V1}. In the last part
of this section we present the main tools used for the proof of
the main theorem.

We write {\it $n$-polytope} instead of "$n$-dimensional polytope" for short.
By {\it facet} we mean a face of codimension one.

\subsection{Gale diagrams}
\label{sec-Gale}

An $n$-polytope is called {\it simple} if any its $k$-face belongs
to exactly $n-k$ facets. Proposition~\ref{3.1} implies that any
compact hyperbolic Coxeter polytope is simple. From now on we
consider simple polytopes only.

Every combinatorial type of simple $n$-polytope with $d$ facets can be
represented by its {\it Gale diagram} $\D$. This consists of $d$ points
$a_1,\dots,a_d$ on the $(d-n-2)$-dimensional unit sphere in
$\R^{d-n-1}$ centered  at the origin.

The combinatorial type of a simple convex polytope can be read off from
the Gale diagram in the following way.
Each point $a_i$ corresponds to the facet
$f_i$ of $P$. For any subset $J$ of the set of facets of $P$ the
intersection of facets $\{f_j \,|\, j\in J \}$ is a face of $P$ if
and only if the origin is contained in the interior of
$\conv\{a_{j} \,|\, j\notin J\}$.

The points $a_1,\dots,a_{d}\in \S^{d-n-2}$ compose a Gale diagram of
some $n$-dimensional polytope $P$ with $d$ facets if and only if
every open half-space $H^+$ in $\R^{d-n-1}$ bounded by a hyperplane
$H$ through the origin contains at least two of the points
$a_1,\dots,a_{d}$.

We should notice that the definition of Gale diagram introduced above 
is "dual" to the standard one (see, for example,~\cite{gale}): usually 
Gale diagram is defined in terms of vertices of polytope instead of facets.
Notice also that the definition above concerns simple polytopes only, and it 
takes simplices out of consideration: usually one  
means the origin of $\R^1$ with multiplicity $n+1$ by the Gale
diagram of an $n$-simplex, however we exclude the origin since we
consider simple polytopes only, and the origin is not contained in
$\D$ for any simple polytope except simplex.

We say that two Gale diagrams $\D$ and $\D'$ are {\it isomorphic} if the
corresponding polytopes are combinatorially equivalent.

If  $d=n+3$ then the Gale diagram of $P$ is two-dimensional, i.e. nodes
$a_i$ of the diagram lie on the unit circle.

A {\it standard Gale diagram} of simple $n$-polytope with $n+3$ facets
consists of vertices $v_1,\dots,v_k$ of regular $k$-gon ($k$ is
odd) in $\R^2$ centered  at the origin which are labeled
according to the following rules:

1) Each label is a positive integer, the sum of labels equals $n+3$.

2) The vertices that lie in any open half-space bounded by a line
 through the origin have labels whose sum is at least two.

Each point $v_i$ with label $\mu_i$ corresponds to $\mu_i$ facets
$f_{i,1},\dots,f_{i,\mu_i}$ of $P$. For any subset $J$ of the set
of facets of $P$ the intersection of facets $\{f_{j,\gamma} \,|\,
(j,\gamma)\in J \}$ is a face of $P$ if and only if the origin is
contained in the interior of $\conv\{v_{j} \,|\, (j,\gamma)\notin
J\}$.

It is easy to check (see, for example,~\cite[Sec. 6.3]{gale}) that any
two-dimensional Gale diagram is isomorphic to some standard
diagram. Two simple $n$-polytopes with $n+3$ facets are
combinatorially equivalent if and only if their standard Gale
diagrams are congruent.

\subsection{Coxeter diagrams}
\label{sec-Cox}

Any Coxeter polytope $P$ can be represented by its Coxeter diagram.

An {\it abstract Coxeter diagram} is a one-dimensional simplicial complex
with weighted edges, where weights are either of the type
$\cos\frac{\pi}{m}$ for some integer $m\ge 3$ or positive real numbers
no less than one. We can suppress the weights but indicate the same information by
labeling the edges of a Coxeter diagram in the following way:

\noindent$\bullet$
if the weight equals   $\cos\frac{\pi}{m}$ then the nodes are
joined by either an $(m-2)$-fold edge or a simple edge labeled by
$m$;

\noindent$\bullet$
if the weight equals one then the nodes are joined by a bold edge;

\noindent$\bullet$
if the weight is greater than one then the nodes are joined by a
dotted edge labeled by its weight.

A {\it subdiagram} of Coxeter diagram is a subcomplex with the same 
as in $\Sigma$. The {\it order}  $|\Sigma|$ is the number of vertices
of the diagram $\Sigma$.

If $\Sigma_1$ and $\Sigma_2$ are sub\-dia\-grams of a
Coxeter diagram $\Sigma$, we denote by $\l \Sigma_1,\Sigma_2\r$ a
sub\-dia\-gram of $\Sigma$ spanned by all nodes of $\Sigma_1$
and $\Sigma_2$. We say that a node of $\Sigma$ {\it attaches} to a 
subdiagram $\Sigma_1\subset\Sigma$ if it is joined with some nodes 
of $\Sigma_1$ by edges of any type.

Let $\Sigma$ be a diagram with $d$ nodes $u_1$,...,$u_d$. Define a
symmetric $d\times d$ matrix  $\Gr(\Sigma)$ in the following way:
$g_{ii}=1$; if two nodes  $u_i$ and $u_j$ are adjacent then
$g_{ij}$ equals negative weight of the edge $u_iu_j$; if two nodes
$u_i$ and $u_j$ are not adjacent then $g_{ij}$ equals zero.

By signature and determinant of diagram $\Sigma$
we mean the signature and the determinant of the matrix
$\Gr(\Sigma)$.

An abstract Coxeter diagram $\Sigma$ is called {\it elliptic} if
the matrix $\Gr(\Sigma)$ is positive definite. 
A Coxeter diagram $\Sigma$ is called {\it parabolic} if
the matrix $\Gr(\Sigma)$ is degenerate, and any subdiagram of
$\Sigma$ is elliptic. Connected elliptic and parabolic diagrams
were classified by Coxeter~\cite{Cox}. We represent the list in
Table~\ref{el-par}.

%
\begin{table}[!h]
\caption{ 
Connected elliptic and parabolic Coxeter diagrams are
listed in left and right columns respectively.
 }
\label{el-par}
\begin{center}
\begin{tabular}{|cc@{\quad}|cc|}
\hline
\raisebox{0pt}{${  A_n}$ $(n\ge 1)$}  & 
\raisebox{0pt}{\epsfig{file=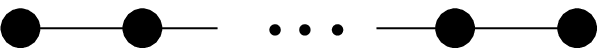,width=0.2\linewidth}}&
\multicolumn{2}{c|}{
\begin{tabular}{cc}
\raisebox{0pt}[20pt][5pt]{${  \widetilde A_1}$} & 
\raisebox{0pt}[20pt][5pt]{\epsfig{file=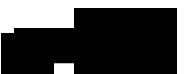,width=0.052\linewidth}}\\
\raisebox{3pt}[15pt][14pt]{${  \widetilde A_n}$ $(n\ge 2)$}  & 
\raisebox{-8pt}[25pt][7pt]{\epsfig{file=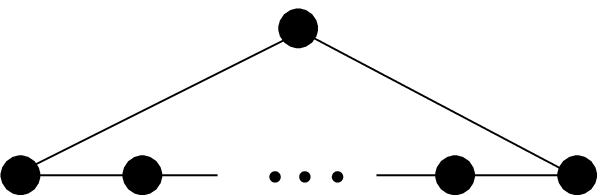,width=0.2\linewidth}}
\end{tabular}
}\\
\hline
\raisebox{-1pt}[23pt][7pt]{${  B_n=C_n}$} & 
\raisebox{-7pt}[23pt][7pt]{\epsfig{file=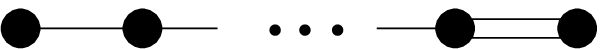,width=0.2\linewidth}}&
\raisebox{7pt}[23pt][7pt]{${  \widetilde B_n}$ $(n\ge 3)$}  & 
\raisebox{-0pt}[30pt][7pt]{\epsfig{file=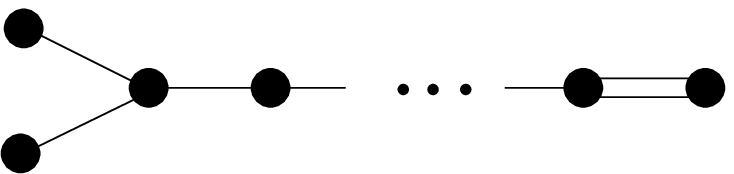,width=0.2\linewidth}}\\
\cline{3-4}
\raisebox{7pt}[15pt][7pt]{ $(n\ge 2)$} & 
&
\raisebox{-0pt}[15pt][7pt]{${  \widetilde C_n}$ $(n\ge 2)$}  & 
\raisebox{-0pt}[15pt][7pt]{\epsfig{file=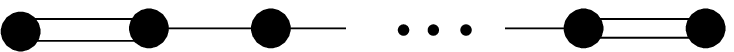,width=0.2\linewidth}}\\
\hline
\raisebox{7pt}[23pt][7pt]{${  D_n}$ $(n\ge 4)$} & 
\raisebox{0pt}[30pt][7pt]{\epsfig{file=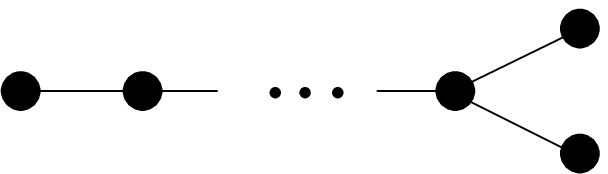,width=0.2\linewidth}}&
\raisebox{7pt}[23pt][7pt]{${  \widetilde D_n}$ $(n\ge 4)$} & 
\raisebox{0pt}[30pt][7pt]{\epsfig{file=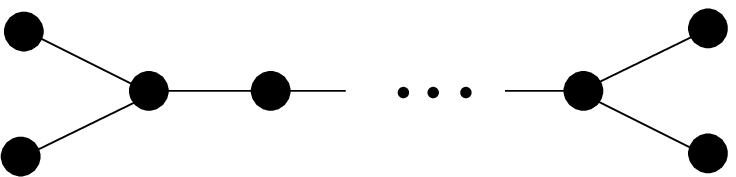,width=0.2\linewidth}}\\
\hline
\raisebox{0pt}[15pt][7pt]{${  G_2^{(m)}}$}  & 
\psfrag{m}{\scriptsize $m$}
\raisebox{0pt}[15pt][7pt]{\epsfig{file=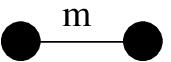,width=0.052\linewidth}}&
\raisebox{0pt}[15pt][7pt]{${  \widetilde G_2}$} & 
\raisebox{0pt}[15pt][7pt]{\epsfig{file=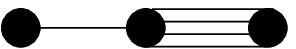,width=0.082\linewidth}}\\
\hline
\raisebox{0pt}[15pt][7pt]{${  F_4}$}  & 
\raisebox{0pt}[15pt][7pt]{\epsfig{file=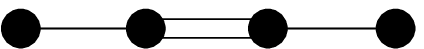,width=0.13\linewidth}}&
\raisebox{-0pt}[15pt][7pt]{${  \widetilde F_4}$} & 
\raisebox{-0pt}[15pt][7pt]{\epsfig{file=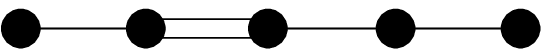,width=0.16\linewidth}}\\
\hline
\raisebox{8pt}[15pt][7pt]{${  E_6}$}  & 
\raisebox{0pt}[30pt][7pt]{\epsfig{file=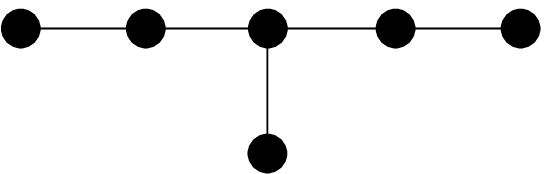,width=0.16\linewidth}}&
\raisebox{8pt}[35pt][7pt]{${  \widetilde E_6}$} & 
\raisebox{-8pt}[40pt][17pt]{\epsfig{file=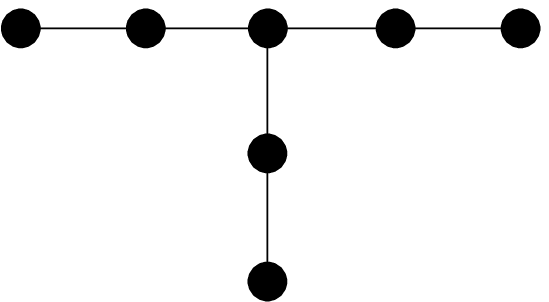,width=0.16\linewidth}}\\
\hline
\raisebox{5pt}[25pt][7pt]{${  E_7}$}  & 
\raisebox{-0pt}[30pt][7pt]{\epsfig{file=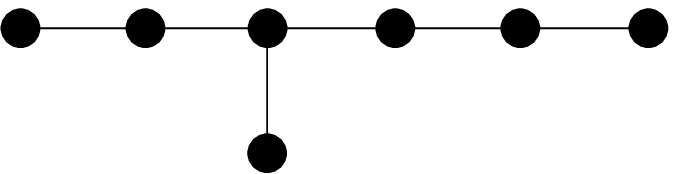,width=0.2\linewidth}}&
\raisebox{5pt}[25pt][7pt]{${  \widetilde E_7}$} & 
\raisebox{-0pt}[25pt][7pt]{\epsfig{file=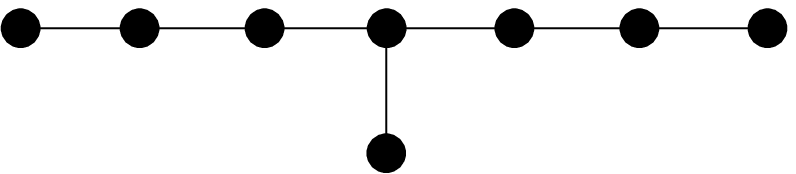,width=0.22\linewidth}}\\
\hline
\raisebox{5pt}[25pt][7pt]{${  E_8}$}  & 
\raisebox{-0pt}[30pt][7pt]{\epsfig{file=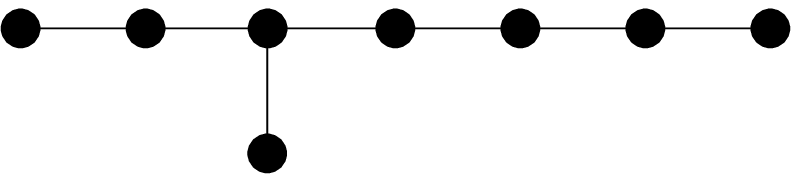,width=0.24\linewidth}}&
\raisebox{5pt}[25pt][7pt]{${  \widetilde E_8}$} & 
\raisebox{-0pt}[25pt][7pt]{\epsfig{file=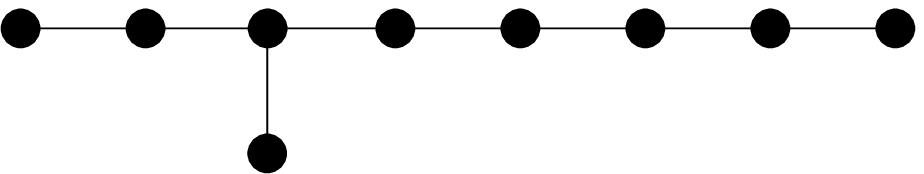,width=0.25\linewidth}}\\
\hline
\raisebox{0pt}[15pt][7pt]{${  H_3}$}  & 
\raisebox{0pt}[15pt][7pt]{\epsfig{file=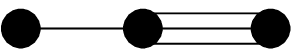,width=0.08\linewidth}}&
& 
\\
\cline{1-2}
\raisebox{0pt}[15pt][7pt]{${  H_4}$}  & 
\raisebox{0pt}[15pt][7pt]{\epsfig{file=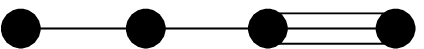,width=0.13\linewidth}}&
& 
\\
\hline
\end{tabular}
\end{center}
\end{table}
%

A Coxeter diagram $\Sigma$ is called a {\it Lann\'er diagram} if
any subdiagram of $\Sigma$ is elliptic, and the diagram $\Sigma$
is neither elliptic nor parabolic. Lann\'er diagrams were
classified by Lann\'er~\cite{L}. We represent the list in
Table~\ref{lan}. A diagram $\Sigma$  is {\it superhyperbolic} if
its negative inertia index is greater than $1$.

%
\begin{table}[!ht]
\caption{{ Lann\'er diagrams.} 
 }
\label{lan}
\begin{center}
\begin{tabular}{|c|cccccc|}
\hline
\raisebox{0pt}{order}  & 
\multicolumn{6}{c|}{
\raisebox{0pt}{diagrams}}\\
\hline
\raisebox{0pt}{\raisebox{-2pt}[10pt][7pt]{$2$}}&
\multicolumn{6}{c|}{\raisebox{0pt}{\epsfig{file=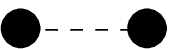,width=0.048\linewidth}}}\\  
\hline
\raisebox{-0pt}{$3$}& 
&\multicolumn{2}{r}{
\psfrag{k}{\scriptsize $k$}
\psfrag{l}{\scriptsize $l$}
\psfrag{m}{\scriptsize $m\;\;\; $}
\begin{tabular}{c}
\raisebox{-0pt}[35pt][10pt]{\epsfig{file=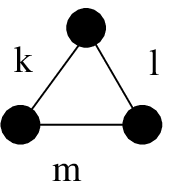,width=0.056\linewidth}}
\end{tabular}
}
&\multicolumn{2}{l}{
\begin{tabular}{l}
$(2\le k,l,m<\infty$,\\
$\frac{1}{k}+\frac{1}{l}+\frac{1}{m}<1)$
\end{tabular}
}&
\\
\hline
\raisebox{-0pt}{$4$}
&\multicolumn{2}{c}{
\begin{tabular}{l}
\raisebox{10pt}{\epsfig{file=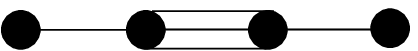,width=0.12\linewidth}}\\
\raisebox{6pt}{\epsfig{file=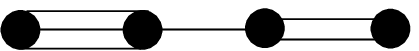,width=0.12\linewidth}}\\
\raisebox{-0pt}{\epsfig{file=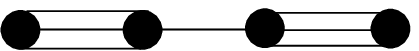,width=0.12\linewidth}}\\
\raisebox{-20pt}{\epsfig{file=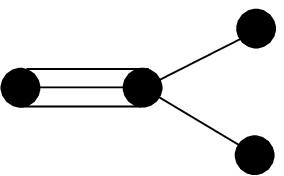,width=0.08\linewidth}}
\end{tabular}
}
&\multicolumn{4}{c|}{
\begin{tabular}{c@{\qquad}c}
\raisebox{-0pt}[30pt][10pt]{\epsfig{file=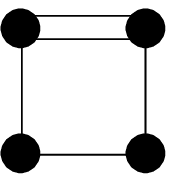,width=0.052\linewidth}}
&\raisebox{-0pt}{\epsfig{file=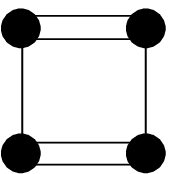,width=0.052\linewidth}}\\
\raisebox{-0pt}[0pt][10pt]{\epsfig{file=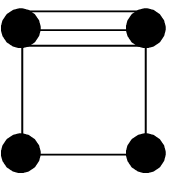,width=0.052\linewidth}}
&\raisebox{-0pt}{\epsfig{file=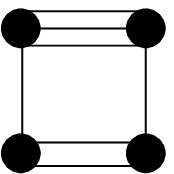,width=0.052\linewidth}}\\
\multicolumn{2}{c}{\raisebox{-0pt}[22pt][10pt]{\epsfig{file=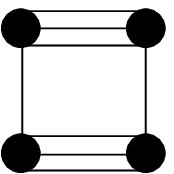,width=0.052\linewidth}}}
\end{tabular}
}\\
\hline
\raisebox{-0pt}{$5$}
&\multicolumn{3}{l}{
\begin{tabular}{l}
\raisebox{7pt}[30pt][0pt]{\epsfig{file=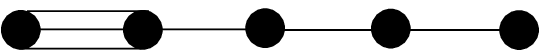,width=0.16\linewidth}}\\
\raisebox{8pt}{\epsfig{file=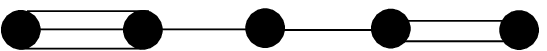,width=0.16\linewidth}}\\
{\raisebox{9pt}{\epsfig{file=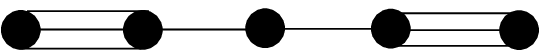,width=0.16\linewidth}}}\\
{\raisebox{8pt}{\epsfig{file=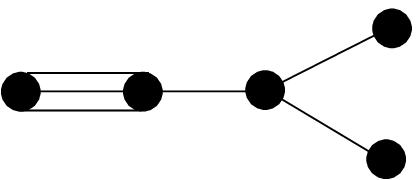,width=0.12\linewidth}}}
\end{tabular}
}
&&\multicolumn{2}{l|}{\raisebox{-0pt}[50pt][10pt]{\epsfig{file=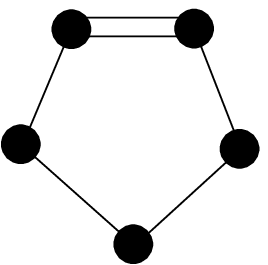,width=0.08\linewidth}}}\\
\hline

\end{tabular}
\end{center}

\end{table}
%

By a {\it simple} (resp., {\it multiple}) edge of Coxeter diagram we mean
an $(m-2)$-fold edge where $m$ is equal to (resp., greater
than) $3$. The number $m-2$ is called the {\it multiplicity} of a
multiple edge. Edges of multiplicity greater than $3$ we call {\it
multi-multiple} edges. If an edge $u_iu_j$ has multiplicity $m-2$
(i.e. the corresponding facets form an angle $\frac{\pi}m$), we
write $[u_i,u_j]=m$. 

A {\it Coxeter diagram $\Sigma(P)$ of Coxeter polytope} $P$ is a Coxeter
diagram whose matrix $\Gr(\Sigma)$ coincides with Gram matrix of
 outer unit normals to the facets of $P$ (referring to the
standard model of hyperbolic $n$-space in $\R^{n,1}$). In other words, 
nodes of Coxeter diagram correspond to facets
of $P$. Two nodes are joined by either an $(m-2)$-fold edge or an
$m$-labeled edge if the corresponding dihedral angle equals
$\frac{\pi}{m}$. If the corresponding facets are parallel the
nodes are joined by a bold edge, and if they diverge then the
nodes are joined by a dotted edge (which may be labeled by
hyperbolic cosine of distance between the hyperplanes containing 
these facets).

If $\Sigma(P)$ is the Coxeter diagram of $P$ then nodes of $\Sigma(P)$ are in
one-to-one correspondence with elements of the set $I=\{1,\dots,d\}$.
For any subset $J\subset I$ denote by $\Sigma(P)_J$ the subdiagram of
$\Sigma(P)$ that consists of nodes corresponding to elements of
$J$. 

\subsection{Hyperbolic Coxeter polytopes}
\label{cox}

In this section by polytope we mean a (probably non-compact) intersection of
closed half-spaces.

\begin{prop}[\cite{V1}, Th. 2.1]
\label{2.1}
Let $\Gr=(g_{ij})$ be indecomposable symmetric matrix of signature
$(n,1)$, where $g_{ii}=1$ and $g_{ij}\le 0$ if  $i\ne j$. Then
there exists a unique (up to isometry of $\H^n$) convex polytope
$P\subset\H^n$ whose Gram matrix coincides with $\Gr$.

\end{prop}

Let $\Gr$ be the Gram matrix of the polytope $P$, and let $J\subset
I$ be a subset of the set of facets of $P$. Denote by $\Gr_J$ the Gram
matrix of vectors $\{e_i\,|\,i\in J\}$, where $e_i$ is outward
unit normal to the facet $f_i$ of $P$ (i.e. $\Gr_J=\Gr(\Sigma(P)_J)$). 
Denote by $|J|$ the number of elements of $J$.

\begin{prop}[\cite{V1}, Th. 3.1]
\label{3.1}
Let $P\subset\H^n$ be an acute-angled polytope with Gram matrix $\Gr$,
and let $J$ be a subset of the set of facets of $P$. The set
$$
q=P\cap\bigcap_{i\in J}f_{i}
$$
is a face of $P$ if and only if the matrix $\Gr_J$ is positive
definite. Dimension of $q$ is equal to $n-|J|$.
\end{prop}

Notice that Prop.~\ref{3.1} implies that the combinatorics of $P$
is completely determined by the Coxeter diagram $\Sigma(P)$.

\smallskip

Let $A$ be a symmetric matrix whose non-diagonal elements are
 non-positive.  $A$ is called
{\it indecomposable} if it cannot be transformed to a
block-diagonal matrix via simultaneous permutations of columns and
rows. We say $A$  to be {\it parabolic} if any indecomposable
 component  of $A$ is positive semidefinite and degenerate. For
 example, a matrix $\Gr(\Sigma)$ for any parabolic diagram
 $\Sigma$ is parabolic.

\begin{prop}[\cite{V1}, cor. of Th. 4.1, Prop. 3.2 and Th. 3.2]
\label{3.2?}
Let $P\subset\H^n$ be a compact Coxeter polytope, and let $\Gr$ be its
Gram matrix. Then for any $J\subset I$ the matrix $\Gr_J$ is not
parabolic.
\end{prop}

Corollary~\ref{cl} reformulates Prop.~\ref{3.2?} in terms of
Coxeter diagrams.

\begin{cor}
\label{cl}
Let $P\subset\H^n$ be a compact Coxeter polytope, and let $\Sigma$ be its
Coxeter matrix. Then any non-elliptic subdiagram of $\Sigma$ contains a
Lann\'er subdiagram.

\end{cor}

\begin{prop}[\cite{V1}, Prop. 4.2]
\label{4.2}
A polytope $P$ in $\H^{n}$ is compact if and only if
it is combinatorially equivalent to some compact convex
$n$-polytope.

\end{prop}

The main result of paper~\cite{nodots} claims that if $P$ is a compact 
hyperbolic Coxeter $n$-polytope having no pair of disjoint facets, then $P$ 
is either a simplex or one of the seven polytopes with $n+2$ facets described
in~\cite{Ess2}. As a corollary, we obtain the following proposition. 

\begin{prop}
\label{nodots}
Let $P\subset\H^n$ be a compact Coxeter polytope with at least
$n+3$ facets. Then $P$ has a pair of disjoint facets.

\end{prop}

\subsection{Coxeter diagrams, Gale diagrams, and missing faces}

Now, for any compact hyperbolic Coxeter polytope we have two diagrams
which carry the complete information about its combinatorics,
namely Gale diagram and Coxeter diagram. The interplay between
them is described by the following lemma, which is a reformulation
of results listed in Section~\ref{cox} in terms of Coxeter
diagrams and Gale diagrams.

\begin{lemma}
\label{face}
A Coxeter diagram $\Sigma$ with nodes $\{u_i\,|\,i=1,\dots,d\}$ is a Coxeter
diagram of some compact hyperbolic Coxeter $n$-polytope with $d$
facets if and only if the following two conditions hold:

1) $\Sigma$ is of signature $(n,1,d-n-1)$;

2) there exists a $(d-n-1)$-dimensional Gale diagram with nodes
   $\{{v_i}\,|\,i=1,\dots,d\}$ and one-to-one map $\psi
   :\{u_i\,|\,i=1,\dots,d\}\to\{{v_i}\,|\,i=1,\dots,d\}$ such that for any
   $J\subset\{1,\dots,d\}$
  the subdiagram $\Sigma_J$ of $\Sigma$ is elliptic if
   and only if the origin is contained in the interior of
$\conv\{\psi({v_i})\,|\,i\notin J\}$.

\end{lemma}

Let $P$ be a simple polytope.
The facets $f_1,\dots,f_m$ of $P$ compose a {\it missing face} of $P$
if $\bigcap\limits_{i=1}^m f_i=\emptyset$ but any proper subset of
$\{f_1,\dots,f_m\}$ has a non-empty intersection.

\begin{prop}[\cite{nodots}, Lemma 2]
\label{missing}
Let $P$ be a simple $d$-polytope with $d+k$ facets $\{f_i\}$, let
$G=\{a_i\}\subset\S^{k-2}$ be a Gale diagram of $P$, and let
$I\subset\{1,\dots,d+k\}$.
Then the set $M_I=\{f_i\,|\,i\in I\}$ is a missing face of
$P$ if and only if the following two conditions hold:
\begin{itemize}
\item[(1)] there exists a hyperplane $H$ through the origin separating
the set $\wh M_I=\{a_i\,|\,i\in I\}$ from the remaining points of $G$;
\item[(2)] for any proper subset $J\subset I$ no hyperplane through
the origin separates the set $\wh M_J=\{a_i\,|\,i\in J\}$ from the
remaining points of $G$.

\end{itemize}
\end{prop}

\noindent
{\bf Remark.}\
Suppose that $P$ is a compact hyperbolic Coxeter polytope.
The definition of missing face (together with Cor.~\ref{cl}) implies
that for any Lann\'er sub\-dia\-gram $L\subset\Sigma(P)$ the
facets corresponding to $L$ compose a missing face of $P$, and any
missing face of $P$ corresponds to some Lann\'er diagram in
$\Sigma(P)$.\\

Now consider a compact hyperbolic Coxeter
$n$-polytope $P$ with $n+3$ facets with standard
Gale diagram $\D$ (which is a $k$-gon, $k$ is odd) and Coxeter diagram
$\Sigma$. Denote by $\Sigma_{i,j}$ a subdiagram of $\Sigma$
corresponding to $j-i+1\pmod k$ consecutive nodes $a_i,\dots,a_j$
of $\D$ (in the sense of Lemma~\ref{face}). If $i=j$, denote
$\Sigma_{i,i}$ by $\Sigma_i$.

The following lemma is an immediate corollary of
Prop.~\ref{missing}.

\begin{lemma}
\label{L}
For any $i\in\{0,\dots,k-1\}$ a diagram $\Sigma_{i+1,i+\frac{k-1}2}$ is a
Lann\'er diagram. All Lann\'er diagrams contained in $\Sigma$ are
of this type.
\end{lemma}

It is easy to see that the collection of missing faces completely
determines the combinatorics of $P$. In view of Lemma~\ref{L} and
the remark above, this means that in Lemma~\ref{face} for given
Coxeter diagram we need to check the signature and correspondence
of Lann\'er diagrams to missing faces of some Gale diagram.\\

\begin{center}
\begin{figure}[!h]
\psfrag{(a)}{\small (a)}
\psfrag{(b)}{\small (b)}
\psfrag{8}{\small $8$}
\psfrag{1}{\small $1$}
\psfrag{2}{\small $2$}
\psfrag{u1}{\small $u_1$}
\psfrag{u2}{\small $u_2$}
\psfrag{u3}{\small $u_3$}
\psfrag{u4}{\small $u_4$}
\psfrag{u5}{\small $u_5$}
\psfrag{u6}{\small $u_6$}
\psfrag{u7}{\small $u_7$}
\begin{center}
\epsfig{file=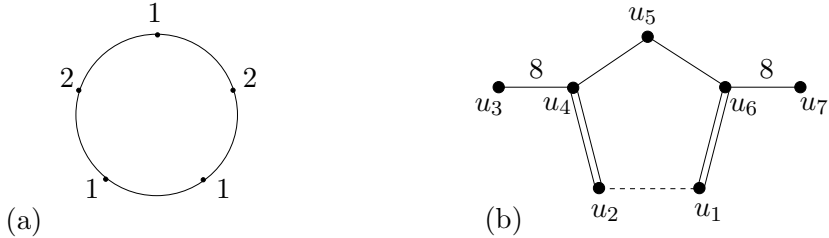,width=0.7\linewidth}
\end{center}
\qquad\qquad
%
\caption{(a) A standard Gale diagram $\D$ and (b) a Coxeter diagram of one of 
polytopes with Gale diagram $\D$}
\label{ex}
\end{figure}
\end{center}

\noindent
{\bf Example.}\ Suppose that there exists a compact hyperbolic
Coxeter polytope $P$ with standard Gale diagram $\D$ shown in
Fig.~\ref{ex}(a). What can we say about Coxeter diagram
$\Sigma=\Sigma(P)$?

The sum of labels of nodes of Gale diagram $\D$ is equal to $7$,
so $P$ is a $4$-polytope with $7$ facets. Thus, $\Sigma$ is
spanned by nodes $u_1,\dots,u_7$, and its signature equals
$(4,1,2)$. Further, $\D$ is a pentagon. By
Lemma~\ref{L}, $\Sigma$ contains exactly $5$ Lann\'er diagrams,
namely $\l u_1,u_2\r$, $\l u_2,u_3,u_4\r$, $\l u_3,u_4,u_5\r$, $\l
u_5,u_6,u_7\r$, and $\l u_6,u_7,u_1\r$.

Now consider the Coxeter diagram $\Sigma$ shown in Fig.~\ref{ex}(b).
Assigning label $1+\sqrt{2}$ to the dotted edge of $\Sigma$, we
obtain a diagram of signature $(4,1,2)$ (this may be shown by
direct calculation). Therefore, there exist $7$ vectors in $\H^4$
with Gram matrix $\Gr(\Sigma)$. It is easy to see that $\Sigma$
contains exactly $5$ Lann\'er diagrams described above. Thus,
$\Sigma$ is a Coxeter diagram of some compact $4$-polytope with
Gale diagram $\D$.

Of course, $\Sigma$ is just an {\it example} of a Coxeter diagram
satisfying both conditions of Lemma~\ref{face} with respect to
given Gale diagram $\D$. In the next two sections we will show how
to list {\it all} compact hyperbolic Coxeter polytopes of given
combinatorial type.

\section{Technical tools}
\label{tech}
From now on by polytope we mean a compact hyperbolic Coxeter
$n$-polytope with $n+3$ facets, and we deal with standard
Gale diagrams only.

\subsection{Admissible Gale diagrams}
\label{gales}
Suppose that there exists a compact hyperbolic Coxeter polytope
$P$ with Gale diagram $\D$. Since the maximal order of
Lann\'er diagram equals five, Lemma~\ref{L} implies that the sum
of labels of $\frac{k-1}2$ consecutive nodes of Gale diagram does
not exceed five. On the other hand, by Lemma~\ref{nodots}, $P$ has
a missing face of order two. This is possible in two cases only:
either $\D$ is a pentagon with two neighboring vertices labeled
by $1$, or $\D$ is a triangle one of whose vertices is labeled by
$2$ (see Prop.~\ref{missing}). Table~\ref{long} contains all
Gale diagrams satisfying one of two conditions above with at least
$7$ and at most $10$ vertices, i.e. Gale diagrams that may
correspond to compact hyperbolic Coxeter $n$-polytopes with $n+3$
facets for $4\le n\le 7$.

\begin{table}[!h]
\caption{Gale diagrams that may correspond to compact Coxeter
polytopes (see Section~\ref{gales})}
\label{long}
\bigskip
$\underline{n=4}$\\

\medskip
\psfrag{1}{\footnotesize $1$}
\psfrag{2}{\footnotesize $2$}
\psfrag{3}{\footnotesize $3$}
\psfrag{4}{\footnotesize $4$}
\psfrag{g1}{\small $\D_1$}
\psfrag{g2}{\small $\D_2$}
\psfrag{g3}{\small $\D_3$}
\psfrag{g4}{\small $\D_4$}
\epsfig{file=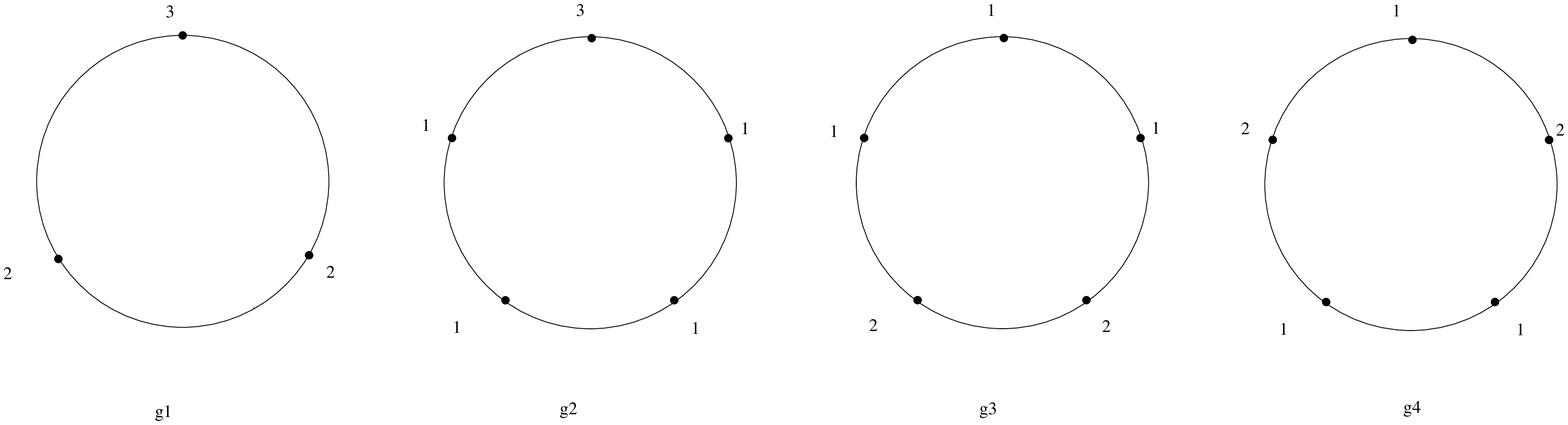,width=0.65\linewidth}\\

\smallskip
$\underline{n=5}$
\begin{center}
\psfrag{1}{\footnotesize $1$}
\psfrag{2}{\footnotesize $2$}
\psfrag{3}{\footnotesize $3$}
\psfrag{4}{\footnotesize $4$}
\psfrag{g5}{\small $\D_5$}
\psfrag{g6}{\small $\D_6$}
\psfrag{g7}{\small $\D_7$}
\psfrag{g8}{\small $\D_8$}
\psfrag{g9}{\small $\D_9$}
\psfrag{g10}{\small $\D_{10}$}
\epsfig{file=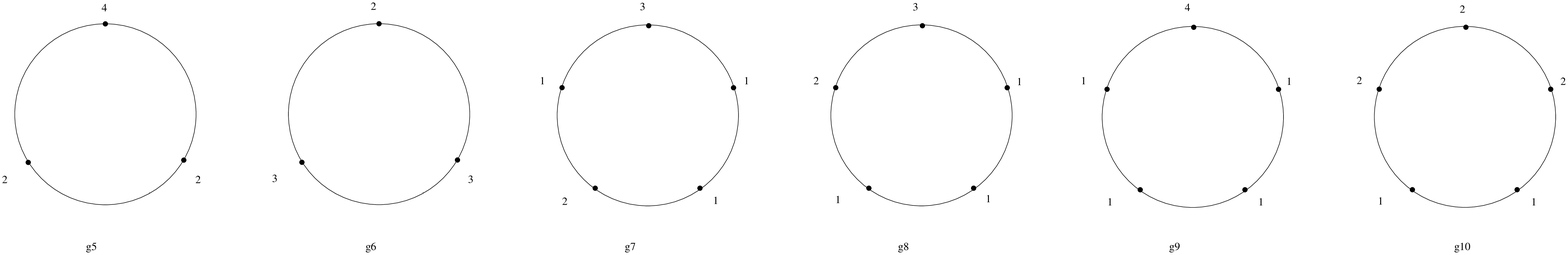,width=0.99\linewidth}
\end{center}
$\underline{n=6}$
\begin{center}
\psfrag{1}{\footnotesize $1$}
\psfrag{2}{\footnotesize $2$}
\psfrag{3}{\footnotesize $3$}
\psfrag{4}{\footnotesize $4$}
\psfrag{5}{\footnotesize $5$}
\psfrag{g11}{\small $\D_{11}$}
\psfrag{g12}{\small $\D_{12}$}
\psfrag{g13}{\small $\D_{13}$}
\psfrag{g14}{\small $\D_{14}$}
\psfrag{g15}{\small $\D_{15}$}
\psfrag{g16}{\small $\D_{16}$}
\epsfig{file=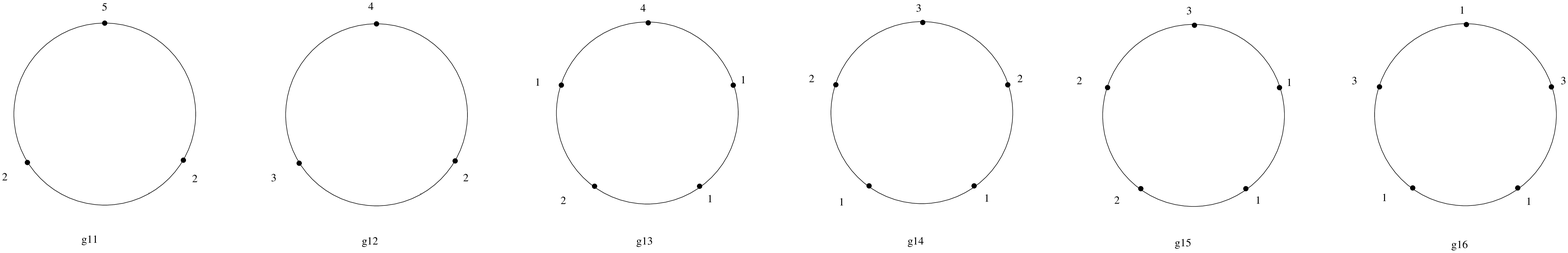,width=0.99\linewidth}
\end{center}
$\underline{n=7}$\\

\smallskip
\psfrag{1}{\footnotesize $1$}
\psfrag{2}{\footnotesize $2$}
\psfrag{3}{\footnotesize $3$}
\psfrag{4}{\footnotesize $4$}
\psfrag{5}{\footnotesize $5$}
\psfrag{g17}{\small $\D_{17}$}
\psfrag{g18}{\small $\D_{18}$}
\psfrag{g19}{\small $\D_{19}$}
\psfrag{g20}{\small $\D_{20}$}
\epsfig{file=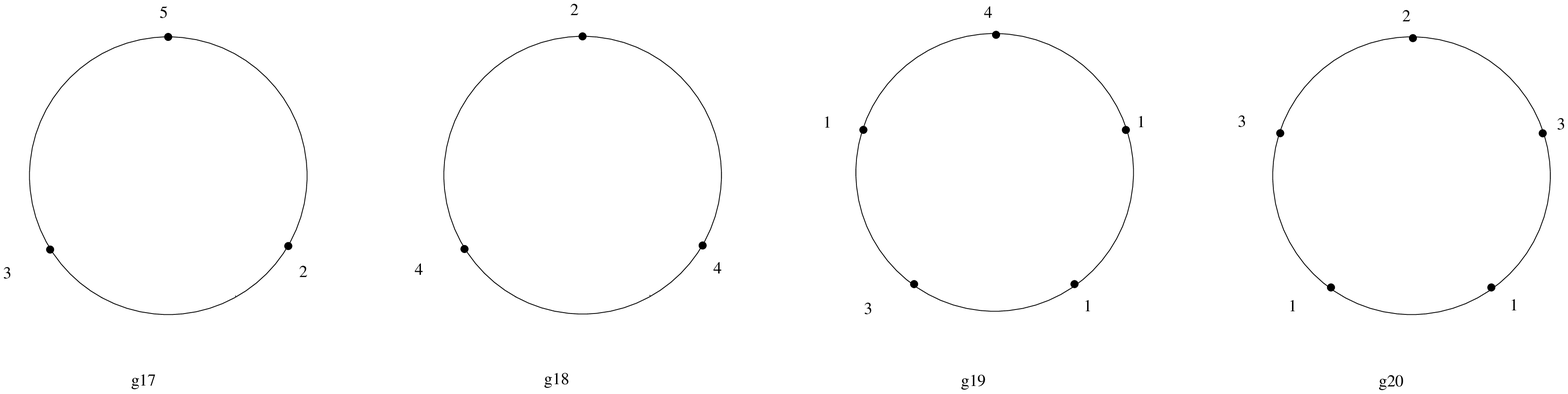,width=0.65\linewidth}
\end{table}

\subsection{Admissible arcs}
\label{arcs}

Let $P$ be an $n$-polytope with $n+3$ facets and let $\D$ be its
$k$-angled Gale diagram. By Lemma~\ref{L}, for any
$i\in\{0,\dots,k-1\}$ the diagram $\Sigma_{i+1,i+\frac{k-1}2}$ is
a Lann\'er diagram. Denote by
$$
\lf x_1,\dots,x_l\rf _{\frac{k-1}2},\quad  l\le k
$$
an arc of length $l$ of $\D$ that consists of $l$ consecutive
nodes with labels $x_1,\dots,x_l$. By writing $J=\lf
x_1,\dots,x_l\rf _{\frac{k-1}2}$ we mean that $J$ is the set of
facets of $P$ corresponding to these nodes of $\D$. The index
$\frac{k-1}2$ means that for any $\frac{k-1}2$ consecutive nodes
of the arc (i.e. for any arc ${I}=\lf
x_{i+1},\dots,x_{i+\frac{k-1}2}\rf _{\frac{k-1}2}$) the subdiagram
$\Sigma_{I}$ of $\Sigma(P)$ corresponding to these nodes is a
Lann\'er diagram (i.e. $I$ is a missing face of $P$).

By Cor.~\ref{cl}, any diagram  $\Sigma_J\subset \Sigma(P)$
corresponding to an arc $J=\lf x_1,\dots,x_l\rf _{\frac{k-1}2}$
satisfies the following property: any subdiagram of $\Sigma_J$
containing no Lann\'er diagram is elliptic. Clearly, any
subdiagram of $\Sigma(P)$ containing at least one
Lann\'er diagram is of signature $(k,1)$ for some $k\le n$.
As it is shown in~\cite{Ess2}, for some arcs $J$ there exist a few
corresponding diagrams $\Sigma_J$ only. In the following lemma, we
recall some results of Esselmann~\cite{Ess2} and prove similar
facts concerning some arcs of Gale diagrams listed in
Table~\ref{long}. This will help us to restrict the number of
Coxeter diagrams that may correspond to some of Gale diagrams listed
in Table~\ref{long}.

\begin{lemma}
\label{dugi}
The diagrams presented in the middle column of Table~\ref{dugi-t} are
the only diagrams that may correspond to arcs listed in the left column.

\end{lemma}


\begin{table}[!h]
\caption{{\normalsize White nodes correspond to endpoints of arcs having
multiplicity one }}
\label{dugi-t}
\smallskip
\begin{center}
\begin{tabular}{|c|c|c|c|}
\hline
&$\ J\ $&\quad all possibilities for $\Sigma_J$\quad &\quad reference (if any)\quad  \\
\hline
1 &
\begin{tabular}{c}
$\lf x,y\rf _1$,\\
$x\ge 4,y\ge 3$ 
\end{tabular}&
\raisebox{0pt}[10pt][8pt]{$\varnothing$}&
\cite{Ess2}, Lemma 4.7\\
\hline
2 &$\lf 1,4,1\rf _2$ &
\begin{tabular}{c}
\raisebox{-3pt}[22pt][5pt]{
\epsfig{file=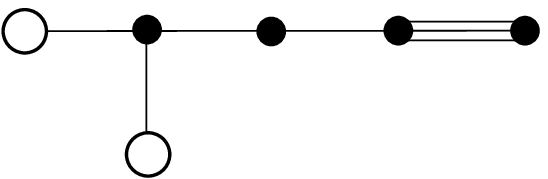,width=0.15\linewidth}
}\\
\epsfig{file=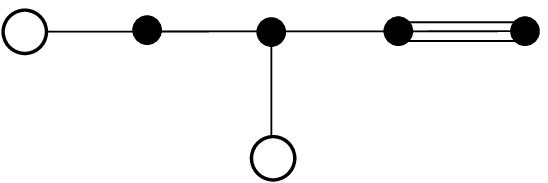,width=0.15\linewidth}\\
\raisebox{7pt}[4pt][0pt]{
\epsfig{file=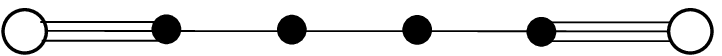,width=0.2\linewidth}
}
\end{tabular}
&
\cite{Ess2}, Lemma 5.3\\
\hline
3 &$\lf 3,2,2\rf _2$&
\raisebox{0pt}[18pt][13pt]{$\varnothing$}
& \cite{Ess2}, Lemma 5.7
\\
\hline
4 &$\lf 4,1,3\rf _2$ &
\raisebox{-7pt}[17pt][13pt]{
\epsfig{file=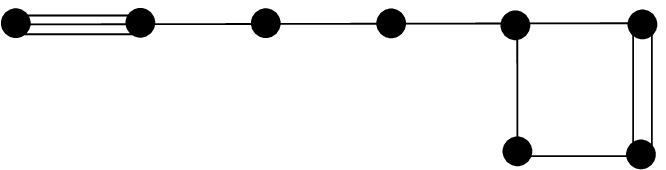,width=0.18\linewidth}
} &
\cite{Ess2}, Lemma 5.9\\
\hline
5 &$\lf 3,1,4,1\rf _2$ &\raisebox{0pt}[11pt][7pt]{$\varnothing$} &
\cite{Ess2}, Folgerung 5.10\\
\hline
6 &$\lf 2,3,2\rf _2$&
\raisebox{0pt}[18pt][13pt]{
\epsfig{file=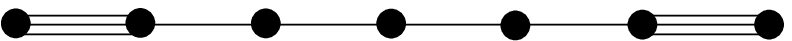,width=0.24\linewidth}
} &
\cite{Ess2}, Lemma 5.12
\\
\hline
7 &$\lf 3,2,3\rf _2$ &
\raisebox{1pt}[11pt][9pt]{
\epsfig{file=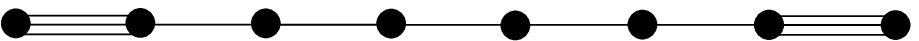,width=0.28\linewidth}
} &
\cite{Ess2}, Lemma 5.12
\\
\hline
$8$ &$\lf 1,3,1\rf _2$ &
\begin{tabular}{c}
\raisebox{1pt}[17pt][7pt]{
\psfrag{3,4}{\scriptsize $3,4$}
\psfrag{4,5}{\scriptsize $4,\!5$}
\psfrag{3,4,5}{\scriptsize $3,4,5$}
\epsfig{file=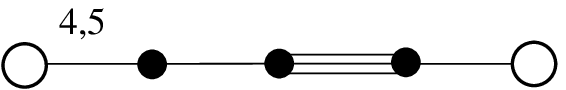,width=0.16\linewidth}
}\\
\begin{tabular}{cc}
\psfrag{3,4}{\scriptsize $3,4$}
\psfrag{4,5}{\scriptsize $4,\!5$}
\psfrag{3,4,5}{\scriptsize $3,4,5$}
\epsfig{file=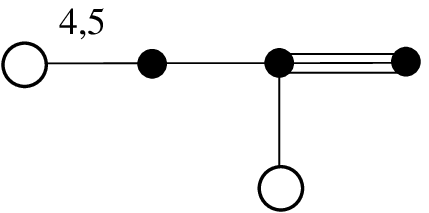,width=0.12\linewidth}
&
\psfrag{3,4}{\scriptsize $3,4$}
\psfrag{4,5}{\scriptsize $4,\!5$}
\psfrag{3,4,5}{\scriptsize $3,4,5$}
\epsfig{file=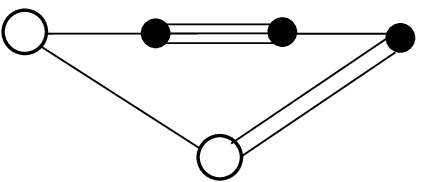,width=0.12\linewidth}
\end{tabular}
\end{tabular}
&\\
\hline
$9$ &
 $\lf 1,3,2\rf _2$
&\begin{tabular}{c}
\raisebox{-5pt}[30pt][5pt]{
\psfrag{3,4}{\scriptsize $3,\!4$}
\psfrag{4,5}{\scriptsize $4,\!5$}
\psfrag{3,4,5}{\scriptsize $3, 4, 5$}
\epsfig{file=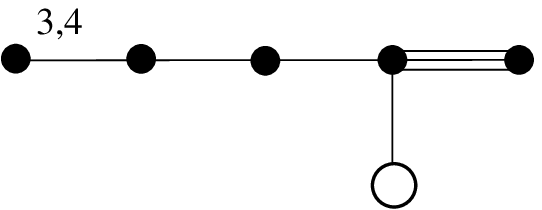,width=0.16\linewidth}
}
\\
\psfrag{3,4}{\scriptsize $3,4$}
\psfrag{4,5}{\scriptsize $4,5$}
\psfrag{3,4,5}{\scriptsize $3,\!4,\!5$}
\epsfig{file=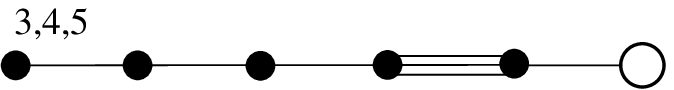,width=0.2\linewidth}\\
\begin{tabular}{cc}
\psfrag{3,4}{\scriptsize $3,\!4$}
\psfrag{4,5}{\scriptsize $4,5$}
\psfrag{3,4,5}{\scriptsize $3,4,5$}
\epsfig{file=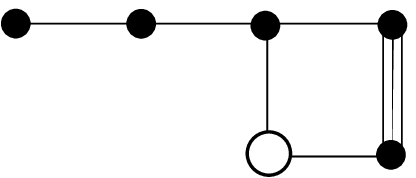,width=0.12\linewidth}&
\epsfig{file=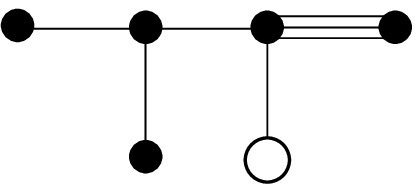,width=0.12\linewidth}
\end{tabular}
\\
\epsfig{file=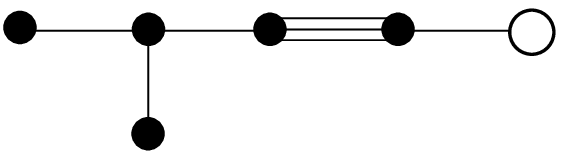,width=0.16\linewidth}
\\
\epsfig{file=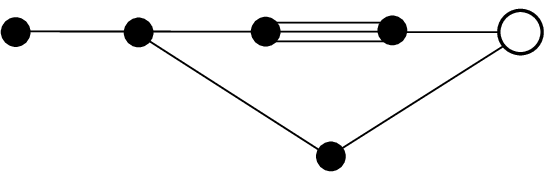,width=0.16\linewidth}
\end{tabular}
&\\
\hline
$10$& $\lf 2,2,2\rf _2$ &
\raisebox{0pt}[18pt][13pt]{
\epsfig{file=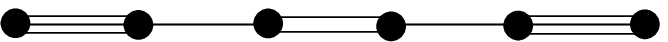,width=0.24\linewidth}
}
&\\
\hline
$11$ &$\lf 3,1,3\rf _2$
&\raisebox{0pt}[10pt][8pt]{$\varnothing$}&\\
\hline
\multicolumn{3}{c}{}
\\
\end{tabular}
\end{center}
\end{table}

\begin{proof}
At first, notice that
for any $J$ as above (i.e. $J$ consists of several consecutive
nodes of Gale diagram) the diagram $\Sigma_J$ must be connected.
This follows from the fact that any Lann\'er diagram is connected,
and that $\Sigma_J$ is not superhyperbolic.

Now we restrict our considerations to items $8$--$11$ only.
For none of these $J$ the diagram $\Sigma_J$ contains a Lann\'er
diagram of order $2$ or $3$. Since $\Sigma_J$ is connected and
does not contain parabolic subdiagrams, this implies that
$\Sigma_J$ does not contain neither dotted nor multi-multiple
edges.
Thus, we are left with finitely many possibilities only, that
allows us to use a computer check: there are several (from $5$ to
$7$) nodes, some of them joined by edges of multiplicity at most
$3$. We only need to check all possible diagrams for the number of
Lann\'er diagrams of all orders and for parabolic subdiagrams.
Namely, in items $8,10$ and $11$ we look for diagrams of order
$5$, $6$ and $7$ containing exactly $2$ Lann\'er subdiagrams of
order $4$ (and containing neither other Lann\'er diagrams nor
parabolic subdiagrams), and in item $9$ we look for diagrams of
order $6$ containing exactly one
Lann\'er subdiagram of order $4$ and exactly one Lann\'er diagram
of order $5$. Notice also that we do not need to check the
signature of obtained diagrams: all them are certainly
non-elliptic, and since any of them contains exactly two
Lann\'er diagrams which have at least one node in common,
by excluding this node we obtain an elliptic diagram.

However, the computation described above is really huge. In what
follows we describe case-by-case how to reduce these computations
to ones taking a few minutes only.

\smallskip
\noindent$\bullet$
{\bf Item 8 ($J=\lf 1,3,1\rf _2$)}. We may consider $\Sigma_J$ as
a
Lann\'er diagram $L$ of order $4$ together with one vertex
attached to $L$ to compose a unique additional Lann\'er diagram
which should be of order $4$, too. There are $9$ possibilities for
$L$ only (Table~\ref{lan}).

\smallskip
\noindent$\bullet$
{\bf Item 9 ($J=\lf 1,3,2\rf _2$)}. The considerations follow the
preceding ones, but we take as $L$ a Lann\'er diagram of order
$5$. Again, there are few possibilities for $L$ only (namely five:
see Table~\ref{lan}).

\smallskip
\noindent$\bullet$
{\bf Item 10 ($J=\lf 2,2,2\rf _2$)}. Again, $\Sigma_J$ contains a
Lann\'er diagram $L$ of order $4$. One of the two remaining nodes of
$\Sigma_J$ must be attached to $L$. Denote this node by $v$. The
diagram $\l L,v\r\subset \Sigma_J$ consists of five nodes and
contains a unique Lann\'er diagram which is of order $4$. All such
diagrams are listed in~\cite[Lemma 3.8]{Ess2} (see the first two
rows of Tabelle $3$, the case $|{\cal N}_F|=1$, $|{\cal L}_F|=4$).
We reproduce this list in Table~\ref{4pl1}.

\begin{table}[!h]
\caption{One of these diagrams should be contained in $\Sigma_J$ for $J=\lf 2,2,2\rf _2$ }
\label{4pl1}
\begin{center}
\psfrag{7}{\small $u_7$}
\epsfig{file=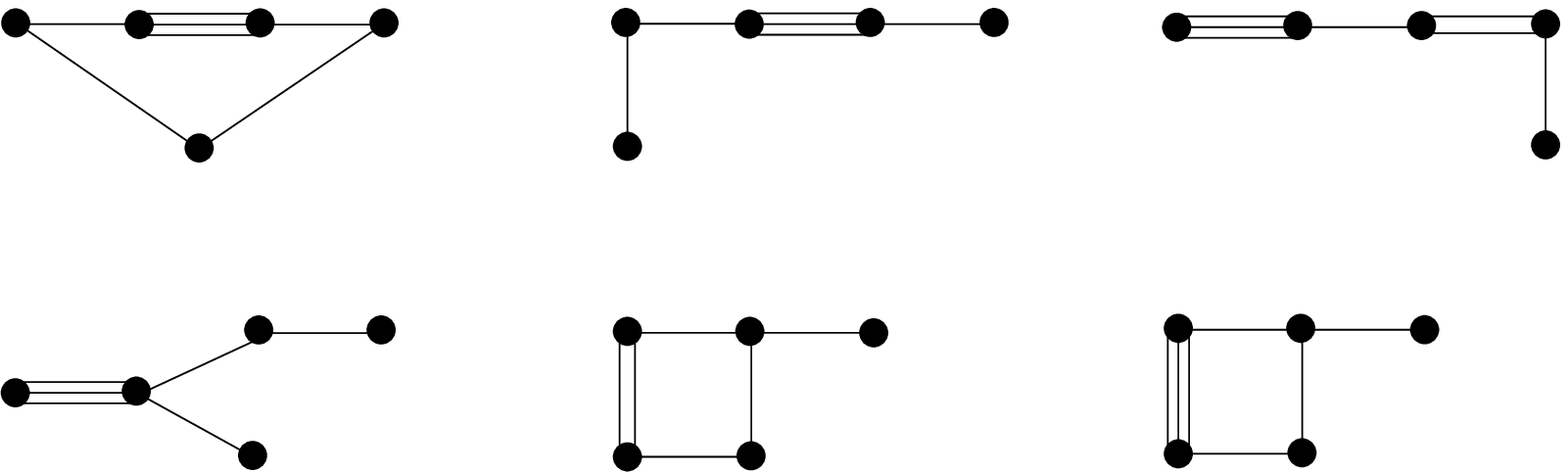,width=0.6\linewidth}
\end{center}
\end{table}

One can see that there are six possibilities only. Now to each of
them we attach the remaining node to compose a unique new Lann\'er
diagram which should be of order $4$.

\smallskip
\noindent$\bullet$
{\bf Item 11 ($J=\lf 3,1,3\rf _2$)}. The considerations are very
similar to the preceding case. $\Sigma_J$ contains a Lann\'er
diagram $L$ of order $4$. One of the three remaining nodes of
$\Sigma_J$ must be attached to $L$. Denote this node by $v$. Now,
one of the two remaining nodes attaches to $\l L,v\r\subset \Sigma_J$.
Denote it by $u$.
The diagram $\l L,v,u\r\subset \Sigma_J$ consists of six nodes and
contains a unique Lann\'er diagram which is of order $4$. All such
diagrams are listed in~\cite[Lemma 3.8]{Ess2} (see Tabelle $3$,
the first two rows of page $27$, the case $|{\cal N}_F|=2$,
$|{\cal L}_F|=4$). We reproduce this list in Table~\ref{4pl2}.

\begin{table}[!h]
\caption{One of these diagrams should be contained in $\Sigma_J$ for $J=\lf 3,1,3\rf _2$}
\label{4pl2}
\begin{center}
\psfrag{7}{\small $u_7$}
\epsfig{file=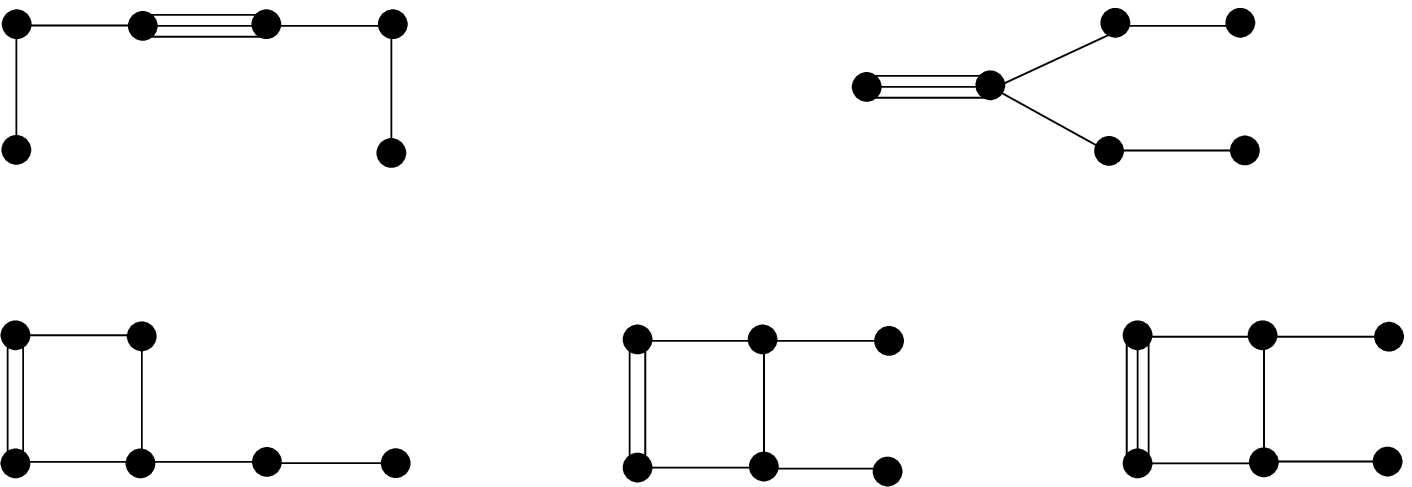,width=0.6\linewidth}
\end{center}
\end{table}

There are five possibilities only. As above, we attach to each of
them the remaining node to compose a unique new Lann\'er diagram which
should be of order $4$.

\end{proof}

\subsection{Local determinants}

In this section we list some tools derived in~\cite{V2} to compute
 determinants of Coxeter diagrams. We will use them to show that
 some (infinite) series of Coxeter diagrams are superhyperbolic.

Let $\Sigma$ be a Coxeter diagram, and let $T$ be a sub\-dia\-gram
of $\Sigma$ such that $\det(\Sigma\setminus T)\ne 0$.
A {\it local determinant} of $\Sigma$ on a sub\-dia\-gram $T$ is
$$\det(\Sigma,T)=\frac{\det \Sigma}{\det(\Sigma\!\setminus\! T)}.$$

\begin{prop}[\cite{V2}, Prop.~12]
\label{loc_sum}
If a Coxeter diagram $\Sigma$ consists of two sub\-dia\-grams
$\Sigma_1$ and $\Sigma_2$ having a unique vertex $v$ in common, 
and no vertex of $\Sigma_1\setminus v$ attaches to
$\Sigma_2\setminus v$, then
$$ \det(\Sigma,v)=\det(\Sigma_1,v)+\det(\Sigma_2,v)-1.$$

\end{prop}

\begin{prop}[\cite{V2}, Prop.~13]
\label{loc_product}
If a Coxeter diagram $\Sigma$ is spanned by two disjoint sub\-dia\-grams
$\Sigma_1$ and $\Sigma_2$ joined by a unique edge $v_1v_2$ of
weight $a$, then
$$\det(\Sigma,\l v_1,v_2\r )=\det(\Sigma_1,v_1)\det(\Sigma_2,v_2) - a^2.$$

\end{prop}

\begin{figure}[!h]
\begin{center}
\psfrag{a}{ $p$}
\psfrag{b}{ $q$}
\psfrag{c}{ $r$}
\psfrag{v}{ $v$}
\epsfig{file=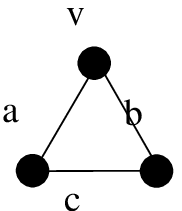,width=0.1\linewidth}
\caption{Diagram $\L_{p,q,r}$}
\label{d_abc}
\end{center}
\end{figure}

Denote by $\L_{p,q,r}$ a Lann\'er diagram of order $3$
containing sub\-dia\-grams of the dihedral groups $G_2^{(p)}$,  $G_2^{(q)}$ and
$G_2^{(r)}$.
Let $v$ be the vertex of $\L_{p,q,r}$ that does not belong to
$G_2^{(r)}$, see Fig.~\ref{d_abc}.
Denote by $\d(p,q,r)$ the local determinant $\det(\L_{p,q,r},v)$.

It is easy to check (see e.g.~\cite{V2}) that
$$
\d(p,q,r)=
1-\frac{\cos^2(\pi/p)+\cos^2(\pi/q)+2\cos(\pi/p)\cos(\pi/q)\cos(\pi/r)}
{\sin^2(\pi/r)}.
$$

Notice that $|\d(p,q,r)|$ is an increasing function on each
of $p,q,r$ tending to infinity while $r$ tends to
infinity.

\section{Proof of the Main Theorem}
\label{proof}

The plan of the proof is the following. First, we show that
there is only a finite number of combinatorial types (or Gale
diagrams) of polytopes we are interested in, and we list these
Gale diagrams. This was done in Table~\ref{long}. For any Gale
diagram from the list we should find all Coxeter polytopes of
given combinatorial type. For that, we try to find all Coxeter
diagrams with the same structure of Lann\'er diagrams as the
structure of missing faces of the Gale diagram is, and then check
the signature. Our task is to be left with finite number of
possibilities for each of Gale diagrams, and use a computer after
that. Some computations involve a large number of cases, but
usually it takes a few minutes of computer's thought. In cases
when it is possible to hugely reduce the computations by better
estimates we do that, but we follow that by long computations to
avoid mistakes.

\begin{lemma}
The following Gale diagrams do not correspond to any
hyperbolic Coxeter polytope: $\D_{12}$, $\D_{15}$, $\D_{16}$,
$\D_{17}$, $\D_{18}$, $\D_{19}$.
\end{lemma}

\begin{proof}
The statement follows from Lemma~\ref{dugi}. Indeed, the diagram
$\D_{12}$ contains an arc $J=\lf 3,4\rf_1$. The corresponding
Coxeter diagram $\Sigma_J$ should be of order $7$, should contain exactly
two Lann\'er diagrams of order $3$ and $4$ which do not intersect,
and should have negative inertia index at most one. Item $1$ of
Table~\ref{dugi-t} implies that there is no such Coxeter diagram
$\Sigma_J$. Thus, $\D_{12}$ is not a Gale diagram of any
hyperbolic
Coxeter polytope.

Similarly, Item $1$ of Table~\ref{dugi-t} also implies the statement of
the lemma for diagrams $\D_{17}$ and $\D_{18}$. Item $3$ implies
the statement for $\D_{15}$, Item $11$ implies the statement for
$\D_{16}$, and Item $5$ implies the statement for the diagram
$\D_{19}$.

\end{proof}

In what follows we check the $14$ remaining Gale diagrams case-by-case.
We start from larger dimensions.

\subsection{Dimension 7}

In dimension $7$ we have only one diagram to consider, namely
$\D_{20}$.

\begin{lemma}
\label{ldim7}
There are no compact hyperbolic Coxeter $7$-polytopes with $10$
facets.

\end{lemma}

\begin{proof}
Suppose that there exists a compact hyperbolic Coxeter polytope $P$ with
Gale diagram $\D_{20}$. This Gale diagram contains an arc $J=\lf 3,2,3\rf_2$.
According to Lemma~\ref{dugi} (Item $7$ of Table~\ref{dugi-t}) and Lemma~\ref{L}, 
the Coxeter diagram $\Sigma$ of $P$ consists of a subdiagram $\Sigma_J$ shown
in Fig.~\ref{323l},
\begin{figure}[!h]
\begin{center}
\psfrag{1}{\small $u_1$}
\psfrag{2}{\small $u_2$}
\psfrag{3}{\small $u_3$}
\psfrag{4}{\small $u_4$}
\psfrag{5}{\small $u_5$}
\psfrag{6}{\small $u_6$}
\psfrag{7}{\small $u_7$}
\psfrag{8}{\small $u_8$}
\epsfig{file=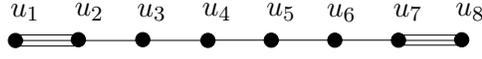,width=0.4\linewidth}
\caption{A unique diagram $\Sigma_J$ for $J=\lf 3,2,3\rf_2$}
\label{323l}
\end{center}
\end{figure}
and two nodes $u_9$, $u_{10}$ joined by a dotted edge. By
Lemma~\ref{face}, the subdiagrams $\l u_{10},u_1,u_2,u_3\r$
and $\l u_6,u_7,u_8,u_9\r$ are Lann\'er diagrams, and no other
Lann\'er subdiagram of $\Sigma$ contains $u_9$ or $u_{10}$.
In particular, $\Sigma$ does not contain Lann\'er subdiagrams of
order $3$.

Consider the diagram $\Sigma'=\l\Sigma_J, u_{9}\r$. It is connected and
contains neither Lann\'er diagrams of order $2$ or $3$, nor
parabolic diagrams. Therefore, $\Sigma'$ does not contain neither
dotted nor multi-multiple edges.
Moreover, by the same reason the node $u_9$ may attach to nodes
$u_1,u_2,u_7$ and $u_8$ by simple edges only. It follows that
there are finitely many possibilities for the diagram $\Sigma'$.
Further, since the diagram $\Sigma'$ defines a collection of $9$ vectors in
$8$-dimensional space $\R^{7,1}$, the determinant of $\Sigma'$ is
equal to zero. A few seconds computer check shows that the only
diagrams satisfying conditions listed in this paragraph are the
following ones:
%
\begin{center}
\psfrag{1}{\small $u_1$}
\psfrag{2}{\small $u_2$}
\psfrag{3}{\small $u_3$}
\psfrag{4}{\small $u_4$}
\psfrag{5}{\small $u_5$}
\psfrag{6}{\small $u_6$}
\psfrag{7}{\small $u_7$}
\psfrag{8}{\small $u_8$}
\psfrag{9}{\small $u_9$}
\epsfig{file=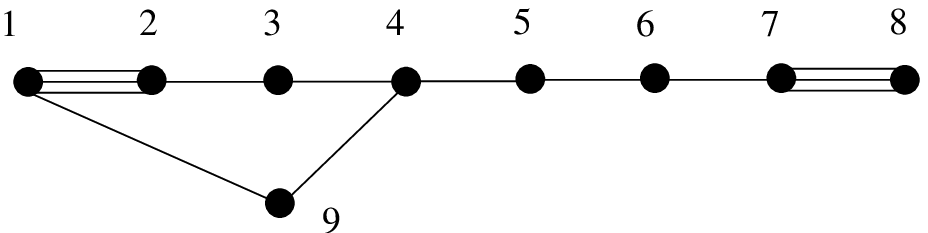,width=0.4\linewidth}
\qquad
\epsfig{file=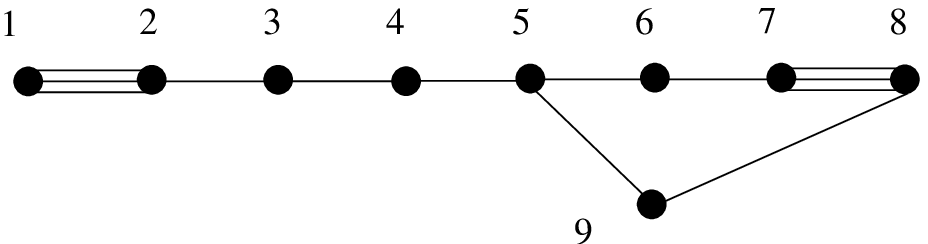,width=0.4\linewidth}
\end{center}
%
However, the left one contains a Lann\'er diagram
$\l u_2,u_1,u_9,u_4,u_5\r$, and the right one contains a Lann\'er
diagram $\l u_7,u_8,u_9,u_5,u_4\r$, which is impossible since $u_9$
does not belong to any Lann\'er diagram of order $5$.

\end{proof}

\subsection{Dimension 6}

In dimension $6$ we are left with three diagrams, namely $\D_{11}$,
$\D_{13}$, and $\D_{14}$.

\begin{lemma}
\label{ldim6-14}
There is only one compact hyperbolic Coxeter polytope
with Gale diagram $\D_{14}$. Its Coxeter diagram is the
lowest one shown in Table~\ref{t6}.

\end{lemma}

\begin{figure}[!h]
\begin{center}
\psfrag{d1}{}
\epsfig{file=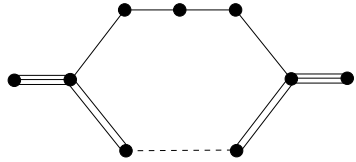,width=0.3\linewidth}
\caption{Coxeter diagram of a unique Coxeter polytope with Gale diagram $\D_{14}$}
\label{6d3s}
\end{center}
\end{figure}

\begin{proof}
Let $P$ be a compact hyperbolic Coxeter polytope with
Gale diagram $\D_{14}$. This Gale diagram contains an arc $J=\lf 2,3,2\rf_2$.
According to Lemma~\ref{dugi} (Item $6$ of Table~\ref{dugi-t}) and Lemma~\ref{L},
the Coxeter diagram $\Sigma$ of $P$ consists of a subdiagram
$\Sigma_J$ shown in Fig.~\ref{232l},
\begin{figure}[!h]
\begin{center}
\psfrag{1}{\small $u_1$}
\psfrag{2}{\small $u_2$}
\psfrag{3}{\small $u_3$}
\psfrag{4}{\small $u_4$}
\psfrag{5}{\small $u_5$}
\psfrag{6}{\small $u_6$}
\psfrag{7}{\small $u_7$}
\epsfig{file=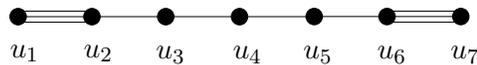,width=0.4\linewidth}
\caption{A unique diagram $\Sigma_J$ for $J=\lf 2,3,2\rf_2$}
\label{232l}
\end{center}
\end{figure}
and two nodes $u_8$, $u_{9}$ joined by a dotted edge.
By Lemma~\ref{face}, the subdiagrams $\l u_{8},u_1,u_2\r$
and $\l u_6,u_7,u_9\r$ are Lann\'er diagrams, and no other
Lann\'er subdiagram of $\Sigma$ contains $u_8$ or $u_9$. So, we
need to check possible multiplicities of edges incident to $u_8$
and $u_9$.

Consider the diagram $\Sigma'=\l\Sigma_J, u_{8}\r$. It is connected,
contains neither Lann\'er diagrams of order $2$ nor parabolic
diagrams, and contains a unique Lann\'er diagram of order $3$,
namely $\l u_{8},u_1,u_2\r$. Therefore, $\Sigma'$ does not contain
dotted edges, and the only multi-multiple edge that may appear
should join $u_8$ and $u_1$.

On the other hand, the signature of $\Sigma_J$ is $(6,1)$. This
implies that the corresponding vectors in $\R^{6,1}$ form a basis,
so the multiplicity of the edge $u_1u_8$ is completely determined
by multiplicities of edges joining $u_8$ with the remaining nodes of
$\Sigma_J$. Since these edges are neither dotted nor
multi-multiple, we are left with a finite number of possibilities
only. We may reduce further computations observing that $u_8$ does
not attach to $\l u_4,u_5,u_6,u_7\r$ (since the diagram $\l u_8,
u_4,u_5,u_6,u_7\r$ should be elliptic), and that multiplicities of
edges $u_8u_2$ and $u_8u_3$ are at most two and one respectively.

Therefore, we have the following possibilities: $[u_8,u_2]=2,3,4$, and,
independently, $[u_8,u_3]=2,3$. For each of these six cases we
should attach the node $u_8$ to $u_1$ satisfying the condition
$\det \Sigma'=0$. An explicit calculation shows that there are two
diagrams listed below.
\begin{center}
\psfrag{1}{\small $u_1$}
\psfrag{2}{\small $u_2$}
\psfrag{3}{\small $u_3$}
\psfrag{4}{\small $u_4$}
\psfrag{5}{\small $u_5$}
\psfrag{6}{\small $u_6$}
\psfrag{7}{\small $u_7$}
\psfrag{8}{\small $u_8$}
\epsfig{file=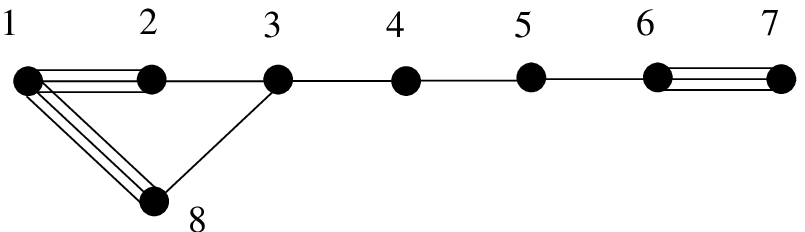,width=0.4\linewidth}
\qquad
\epsfig{file=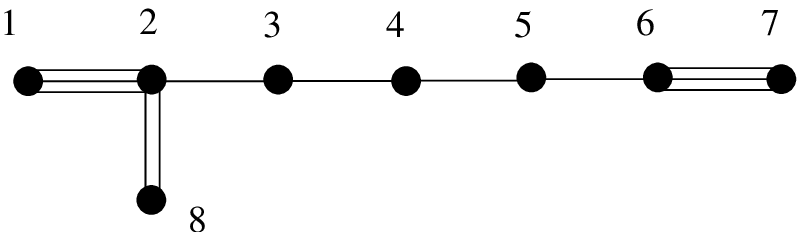,width=0.4\linewidth}
\end{center}
The left one contains a Lann\'er diagram $\l
u_1,u_8,u_3,u_4,u_5\r$, which is impossible. At the same time, the
right one contains exactly Lann\'er diagrams prescribed by Gale
diagram.

Similarly, the node $u_9$ may be attached to $\Sigma_J$ in a
unique way, i.e. by a unique edge $u_9u_6$ of multiplicity two.
Thus, $\Sigma$ must look like the diagram shown in
Fig.~\ref{6d3s}.

Now we write down the determinant of $\Sigma$ as a quadratic polynomial of
the weight $d$ of the dotted edge. An easy computation shows that
$$\det\Sigma=\frac{\sqrt{5}-2}{32}\left(d-(\sqrt{5}+2)\right)^2$$
The signature of $\Sigma$ for $d=\sqrt{5}+2$ is equal to $(6,1,2)$, so we
obtain that this diagram corresponds to a Coxeter polytope.

\end{proof}

\begin{lemma}
\label{ldim6-13}
There are two compact hyperbolic Coxeter polytopes
with Gale diagram $\D_{13}$. Their Coxeter diagrams are
shown in the upper row of Table~\ref{t6}.

\end{lemma}

\begin{proof}
Let $P$ be a compact hyperbolic Coxeter polytope with
Gale diagram $\D_{13}$. This Gale diagram contains an arc $J=\lf 1,4,1\rf_2$.
Hence, the Coxeter diagram $\Sigma$ of $P$ contains a diagram $\Sigma_J$
which coincides with one of the three diagrams shown in Item $2$ of
Table~\ref{dugi-t}. Further, $\Sigma$ contains two Lann\'er diagrams
 of order $3$, one of which (say, $L$) intersects $\Sigma_J$. Denote the
 common node of that Lann\'er diagram $L$ and $\Sigma_J$ by $u_1$,
 the $5$ remaining nodes of $\Sigma_J$ by $u_2,\dots,u_6$ (in a way
 that $u_6$ is marked white in Table~\ref{dugi-t}, i.e. it belongs to
 only one Lann\'er diagram of order $5$), and denote the two remaining 
 nodes of $L$ by $u_7$ and $u_8$. Since $L$ is connected, we may
 assume that $u_7$ is joined with $u_1$. Notice that $u_1$ is also
 a node marked white in Table~\ref{dugi-t}, elsewhere it belongs
 to at least three Lann\'er diagrams in $\Sigma$.

Consider the diagram $\Sigma'=\l\Sigma_J, u_{7}\r$. It is
connected, and all Lann\'er diagrams contained in $\Sigma'$ are
contained in $\Sigma_J$. In particular, $\Sigma'$ does not contain
neither dotted nor multi-multiple edges. Hence, we have only finite
number of possibilities for $\Sigma'$. More precisely, to each of the 
three diagrams $\Sigma_J$ shown in Item $2$ of Table~\ref{dugi-t} we must
attach a node $u_7$ without making new Lann\'er (or parabolic)
diagrams, and all edges must have multiplicities at most $3$. In
addition, $u_7$ is joined with $u_1$. The last condition is
restrictive, since we know that $u_1$ and $u_6$ are the nodes of
$\Sigma_J$ marked white in Table~\ref{dugi-t}. A direct
computation (using the technique described in Section~\ref{arcs})
leads us to the two diagrams $\Sigma'_1$ and $\Sigma'_2$ (up to
permutation of indices $2,3,4$ and $5$ which does not play any
role) shown in Fig.~\ref{1141}.
\begin{figure}[!h]
\begin{center}
\psfrag{1}{\small $u_1$}
\psfrag{2}{\small $u_2$}
\psfrag{3}{\small $u_3$}
\psfrag{4}{\small $u_4$}
\psfrag{5}{\small $u_5$}
\psfrag{6}{\small $u_6$}
\psfrag{7}{\small $u_7$}
\psfrag{8}{\small $u_8$}
$\Sigma'_1=$ \raisebox{-15pt}{\epsfig{file=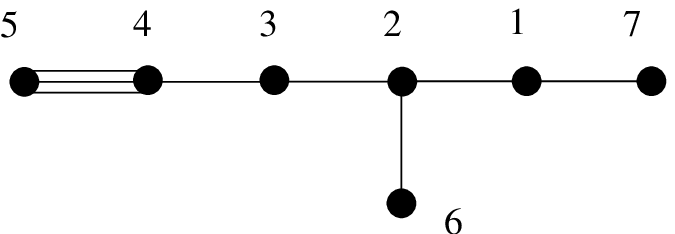,width=0.3\linewidth}}
\qquad\qquad
$\Sigma'_2=$ \raisebox{-15pt}{\epsfig{file=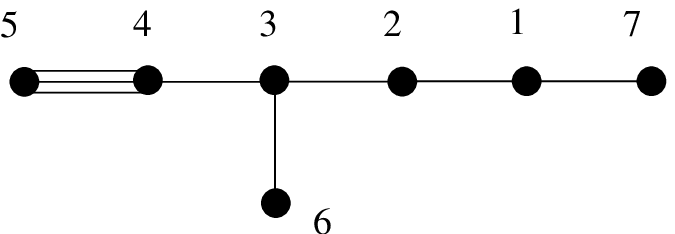,width=0.3\linewidth}}
\caption{Two possibilities for diagram $\Sigma'$, see Lemma~\ref{ldim6-13}}
\label{1141}
\end{center}
\end{figure}

Now consider the diagram $\Sigma''=\l\Sigma', u_{8}\r=\l\Sigma_J, u_{7},
u_8\r=\l\Sigma_J,L\r$. As above, $u_8$ may attach to $\Sigma_J$ by
edges of multiplicity at most $3$, so the only multi-multiple edge
that may appear in $\Sigma''$ is $u_8u_7$. Since both diagrams
$\Sigma'_1$ and $\Sigma'_2$ have signature $(6,1)$, the
corresponding vectors in $\R^{6,1}$ form a basis, so the
multiplicity of the edge $u_8u_7$ is completely determined by
multiplicities of edges joining $u_8$ with the remaining nodes of
$\Sigma'$. Thus, there is a finite number of possibilities for
$\Sigma''$. To reduce the computations note that $u_8$ is not
joined with $\l u_2,u_3,u_4,u_5\r$ (since the diagram $\l
u_2,u_3,u_4,u_5,u_8\r$ must be elliptic). Attaching $u_8$ to
$\Sigma'_2$, we do not obtain any diagram with zero determinant
and prescribed
Lann\'er diagrams. Attaching $u_8$ to $\Sigma'_1$, we obtain the two
diagrams $\Sigma''_1$ and $\Sigma''_2$ shown in Fig.~\ref{2141}.
\begin{figure}[!h]
\begin{center}
\psfrag{1}{\small $u_1$}
\psfrag{2}{\small $u_2$}
\psfrag{3}{\small $u_3$}
\psfrag{4}{\small $u_4$}
\psfrag{5}{\small $u_5$}
\psfrag{6}{\small $u_6$}
\psfrag{7}{\small $u_7$}
\psfrag{8}{\small $u_8$}
\psfrag{10}{\scriptsize $10$}
$\Sigma''_1=$ \raisebox{-15pt}{\epsfig{file=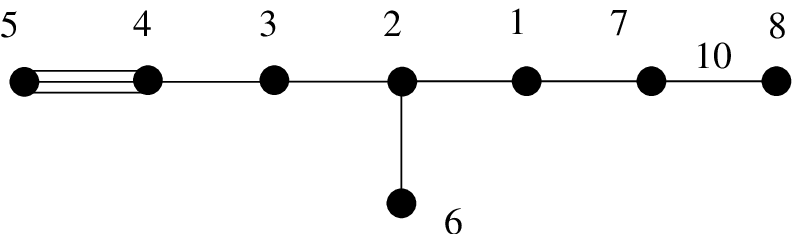,width=0.3\linewidth}}
\qquad\qquad
$\Sigma''_2=$ \raisebox{-15pt}{\epsfig{file=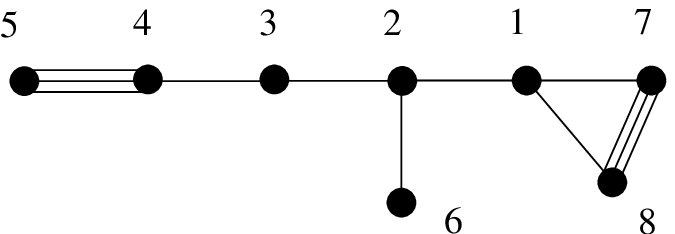,width=0.3\linewidth}}
\caption{Two possibilities for diagram $\Sigma''$, see Lemma~\ref{ldim6-13}}
\label{2141}
\end{center}
\end{figure}

The remaining node of $\Sigma$, namely $u_9$, is joined with $u_6$ by a dotted
edge. It is also contained in a Lann\'er diagram $\l
u_7,u_8,u_9\r$ of order $3$, but no other Lann\'er diagram
contains $u_9$. Since $u_7$ attaches to $u_1$, we see that all
edges joining $u_9$ with $\Sigma'\setminus u_6$ are neither dotted
nor multi-multiple. On the other hand, for both diagrams
$\Sigma''_1$ and $\Sigma''_2$, the diagram $\Sigma''\setminus u_6$
has signature $(6,1)$. Hence, the weight of edge $u_9u_8$ is
completely determined by multiplicities of edges joining $u_9$
with the remaining nodes of $\Sigma''\setminus u_6$, so we are left
with finitely many possibilities for $\Sigma''\setminus u_6$.
Again, we note that $u_9$ is not joined with $\l
u_2,u_3,u_4,u_5\r$. Now we attach $u_9$ to $u_1$ and to $u_7$ by
edges of multiplicities from $0$ (i.e. no edge) to $3$, and then
compute the weight of the edge $u_9u_8$ to obtain
$\det(\Sigma\setminus u_6)=0$. This weight is equal to
$\cos\frac{\pi}m$ for integer $m$ only in case of the diagrams shown
in Fig.~\ref{12141}.
\begin{figure}[!h]
\begin{center}
\psfrag{1}{\small $u_1$}
\psfrag{2}{\small $u_2$}
\psfrag{3}{\small $u_3$}
\psfrag{4}{\small $u_4$}
\psfrag{5}{\small $u_5$}
\psfrag{6}{\small $u_6$}
\psfrag{7}{\small $u_7$}
\psfrag{8}{\small $u_8$}
\psfrag{9}{\small $u_9$}
\psfrag{10}{\scriptsize $10$}
\epsfig{file=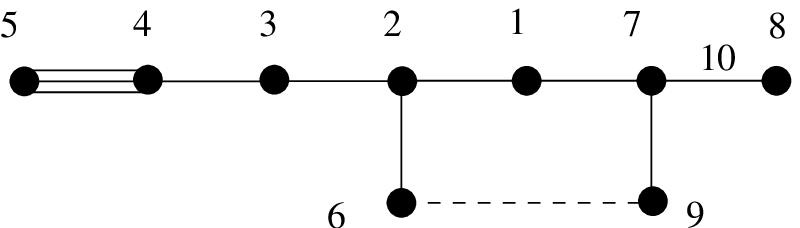,width=0.35\linewidth}
\qquad\qquad
\raisebox{-10pt}{\epsfig{file=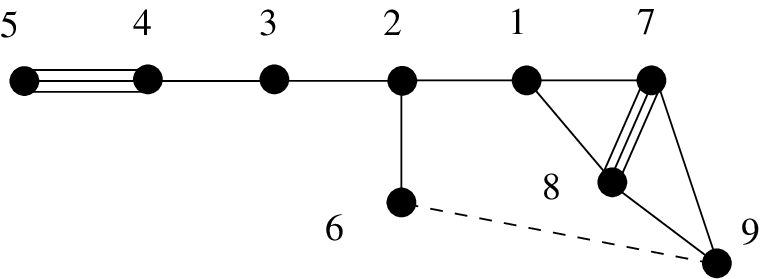,width=0.35\linewidth}}
\caption{Coxeter diagrams of Coxeter polytopes with Gale diagram $\D_{13}$}
\label{12141}
\end{center}
\end{figure}

The last step is to find the weight of the dotted edge $u_9u_6$ to
satisfy the signature condition, i.e. the signature should equal
$(6,1,2)$. We write the determinant of $\Sigma$ as a quadratic
polynomial of the weight $d$ of the dotted edge, and compute the
root. An easy computation shows that for both diagrams the
signature of $\Sigma$ for $d=\frac{1+\sqrt{5}}2$ is equal to
$(6,1,2)$, so we obtain that these two diagrams correspond to
Coxeter polytopes. One can note that the right polytope can be
obtained by gluing two copies of the left one along the facet
corresponding to the node $u_8$.

\end{proof}

\begin{lemma}
\label{ldim6-11}
There are no compact hyperbolic Coxeter polytopes
with Gale diagram $\D_{11}$.

\end{lemma}

\begin{proof}

Suppose that there exists a hyperbolic Coxeter polytope $P$ with
Gale diagram $\D_{11}$. The Coxeter diagram $\Sigma$ of $P$
contains a Lann\'er diagram $L_1=\l u_1,\dots u_5\r$ of order $5$,
and two diagrams of order $2$, denote them $L_2=\l u_6,u_8\r$ and
$L_3=\l u_7,u_9\r$. The diagram $\l L_1,L_2\r$ is connected,
otherwise it is superhyperbolic. Thus, we may assume that $u_6$
attaches to $L_1$. Similarly, we may assume that $u_7$ attaches to
$L_1$.

Therefore, the diagram $\Sigma'=\l L_1,u_6,u_7\r$ consists of a
Lann\'er diagram $L_1$ of order $5$ and two additional nodes which
attach to $L_1$, and these nodes are not contained in any Lann\'er
diagram. According to~\cite[Lemma 3.8]{Ess2} (see Tabelle $3$,
page $27$, the case $|{\cal N}_F|=2$, $|{\cal L}_F|=5$), $\Sigma'$
must coincide with the diagram (up to permutation of indices
of nodes of $L_1$) shown in Fig.~\ref{151}.
\begin{figure}[!h]
\begin{center}
\psfrag{1}{\small $u_1$}
\psfrag{2}{\small $u_2$}
\psfrag{3}{\small $u_3$}
\psfrag{4}{\small $u_4$}
\psfrag{5}{\small $u_5$}
\psfrag{6}{\small $u_6$}
\psfrag{7}{\small $u_7$}
\epsfig{file=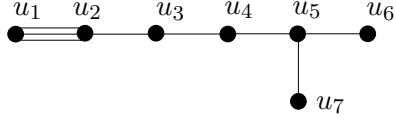,width=0.32\linewidth}
\caption{The diagram $\Sigma'$, see Lemma~\ref{ldim6-11}}
\label{151}
\end{center}
\end{figure}

Consider the diagram $\Sigma''_1=\l\Sigma',u_8\r=\Sigma\setminus
u_9$. The node $u_8$ is joined with $u_6$ by a dotted edge. The
diagram $\Sigma''_1\setminus u_6$ contains a unique Lann\'er
diagram, $L_1$. If $u_8$ attaches to $L_1$, $\Sigma''_1\setminus
u_6$ should coincide with $\Sigma'$. Thus, $u_8$ does not attach
to $\l u_1,\dots,u_4\r$, and $[u_8,u_5]=2$ or $3$. It is also easy
to see that $[u_8,u_7]\le 4$. Since the signature of $\Sigma'$ is
$(6,1)$, the weight of the edge $u_8u_6$ is completely determined
by multiplicities of edges joining $u_8$ with the remaining nodes of
$\Sigma'$. Hence, we have a finite number of possibilities for
$\Sigma''_1$. To reduce the computations observe that either
$[u_5,u_8]$ or $[u_7,u_8]$ must equal $2$. We are left with only
$4$ cases: the pair $([u_5,u_8],[u_7,u_8])$ coincides with one of
$(2,2),(2,3),(2,4)$ or $(3,2)$. For each of them we compute the
weight of $u_8u_6$ by solving the equation $\det\Sigma''_1=0$.
Each of these equations has one positive and one negative
solution, but the positive solution in case of
$([u_5,u_8],[u_7,u_8])=(2,4)$ is less than one, so it cannot be a
weight of a dotted edge. Therefore, we have three cases
$([u_5,u_8],[u_7,u_8])=(2,2),(2,3)$ or $(3,2)$, for which the
weight of $u_8u_6$ is equal to
$\frac{\sqrt{2}\sqrt{4+\sqrt{5}}}{\sqrt{11}}$,
$\frac{-3\sqrt{5}+7+4\sqrt{10-4\sqrt{5}}}{\sqrt{-9+5\sqrt{5}}}$,
 and $\frac{5+4\sqrt{5}}{11}$ respectively.

By symmetry, we obtain the same cases for the diagram
$\Sigma''_2=\l\Sigma',u_9\r=\Sigma\setminus u_8$, and the same
values of the weight of the edge $u_9u_7$ when
$([u_5,u_9],[u_6,u_9])=(2,2),(2,3)$ and $(3,2)$ respectively.
Now, we have only $9$ cases to attach nodes $u_8$ and $u_9$ to
$\Sigma'$ (in fact, there are only six up to symmetry). For each
of these cases we compute the weight of the edge $u_8u_9$ by
solving the equation $\det\Sigma=0$. None of these solutions is
equal to $\cos\frac{\pi}m$ for integer $m$, which contradicts 
the fact that the diagram $\l u_8,u_9\r$ is elliptic. This
contradiction proves the lemma.

\end{proof}

\subsection{Dimension 5}

In dimension $5$ we must consider six Gale diagrams, namely $\D_{5}$-- $\D_{10}$.

\begin{lemma}
\label{ldim5-10}
There is only one compact hyperbolic Coxeter polytope
with Gale diagram $\D_{10}$. Its Coxeter diagram is the
left one shown in the first row of Table~\ref{t5}.

\end{lemma}

\begin{proof}
The proof is similar to the proof of Lemma~\ref{ldim6-14}.
We assume that there exists a hyperbolic Coxeter polytope $P$ with
Gale diagram $\D_{10}$. This Gale diagram contains an arc $J=\lf 2,2,2\rf_2$.
According to Lemma~\ref{dugi} (Item $10$ of Table~\ref{dugi-t}) and Lemma~\ref{L},
the Coxeter diagram $\Sigma$ of $P$ consists of the subdiagram
$\Sigma_J$ shown in Fig.~\ref{222},
\begin{figure}[!h]
\begin{center}
\psfrag{1}{\small $u_1$}
\psfrag{2}{\small $u_2$}
\psfrag{3}{\small $u_3$}
\psfrag{4}{\small $u_4$}
\psfrag{5}{\small $u_5$}
\psfrag{6}{\small $u_6$}
\epsfig{file=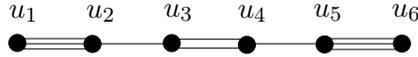,width=0.35\linewidth}
\caption{A unique diagram $\Sigma_J$ for $J=\lf 2,2,2\rf_2$}
\label{222}
\end{center}
\end{figure}
and two nodes $u_7$, $u_{8}$ joined by a dotted edge.
By Lemma~\ref{face}, the subdiagrams $\l u_{7},u_1,u_2\r$
and $\l u_5,u_6,u_8\r$ are Lann\'er diagrams, and no other Lann\'er
subdiagram of $\Sigma$ contains $u_7$ or $u_8$. So, we need to
check possible multiplicities of edges incident to $u_7$ and
$u_8$.

Again, we consider the diagram $\Sigma'=\l\Sigma_J, u_{7}\r$.
It is connected, does not contain dotted edges, and its determinant
is equal to zero. Furthermore, observe that $u_7$ does not attach
to $\l u_2,u_3,u_4,u_5\r$ (since the diagram $\l u_7,
u_2,u_3,u_4,u_5\r$ should be elliptic), and $u_7$ does not attach
to $u_6$ (since the diagram $\l u_7,u_4,u_5,u_6\r$ should be
elliptic). Therefore, $u_7$ is joined with $u_1$ only. Solving the
equation $\det\Sigma'=0$, we find that $[u_7,u_1]=4$.

By symmetry, we obtain that $u_8$ is not joined with $\l u_1,
u_2,u_3,u_4,u_5\r$, and $[u_8,u_6]=4$. Thus, we have the Coxeter
diagram $\Sigma$ shown in Fig.~\ref{12221}.
\begin{figure}[!h]
\begin{center}
\psfrag{d1}{\small $u_2$}
\psfrag{d2}{\small $u_5$}
\psfrag{b1}{\small $u_7$}
\psfrag{b2}{\small $u_8$}
\psfrag{c1}{\small $u_1$}
\psfrag{c2}{\small $u_6$}
\psfrag{e1}{\small $u_3$}
\psfrag{e2}{\small $u_4$}
\epsfig{file=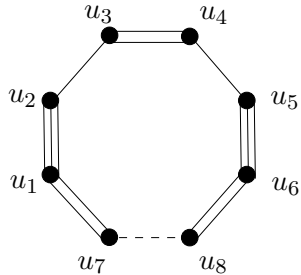,width=0.25\linewidth}
\caption{Coxeter diagram of a unique Coxeter polytope with Gale diagram $\D_{10}$}
\label{12221}
\end{center}
\end{figure}
Assigning the weight $d=\sqrt{2}(\sqrt{5}+1)/4$ to the dotted
edge, we see that the signature of $\Sigma$ is equal to $(5,1,2)$,
so we obtain that this diagram corresponds to a Coxeter polytope.

\end{proof}

Before considering the diagram $\D_{9}$, we make a small geometric
excursus, the first one in this purely geometric paper.

The combinatorial type of polytope defined by Gale diagram
$\D_{9}$ is twice truncated $5$-simplex, i.e. a $5$-simplex in which
two vertices are truncated by hyperplanes very close to the
vertices. If we have such a polytope $P$ with acute angles, it is
easy to see that we are always able to truncate the polytope again
by two hyperplanes in the following way: we obtain a
combinatorially equivalent polytope $P'$; the two truncating
hyperplanes do not intersect initial truncating hyperplanes and
intersect exactly the same facets of $P$ the initial ones do; the
two truncating hyperplanes are orthogonal to all facets of $P$
they do intersect.

The difference between polytopes $P$ and $P'$ consists of two small
polytopes, each of them is combinatorially equivalent to a product
of $4$-simplex and segment, i.e. each of these polytopes is a
simplicial prism. Of course, it is a Coxeter prism, and one of the
bases is orthogonal to all facets of the prism it does intersect.
All such prisms were classified by Kaplinskaja in~\cite{Kap}.
Simplices truncated several times with orthogonality condition
described above were classified by Schlettwein in~\cite{trun}.
Twice truncated simplices from the second list are the right ones
in rows $1$, $3$, and $5$ of Table~\ref{t5}.

Therefore, to classify all Coxeter polytopes with Gale diagram
$\D_9$ we only need to do the following. We take a twice truncated
simplex from the second list, it has two "right" facets, i.e.
facets which make only right angles with other facets. Then we
find all the prisms that have "right" base congruent to one of
"right" facets of the truncated simplex, and glue these prisms to
the truncated simplex by "right" facets in all possible ways.

The result is presented in Table~\ref{t5}. All polytopes except the
left one from the first row have Gale diagram $\D_9$. The
polytopes from the fifth row are obtained by gluing one prism to
the right polytope from this row, the polytopes from the third and
fourth rows are obtained by gluing prisms to the right polytope
from the third row, and the polytopes from the first and second
rows are obtained by gluing prisms to the right polytope from the
first row. The number of glued prisms is equal to the number of
edges inside the maximal cycle of Coxeter diagram. Hence, we come
to the following lemma:

\begin{lemma}
\label{ldim5-9}
There are $15$ compact hyperbolic Coxeter $5$-polytopes with $8$
facets with Gale diagram $\D_{9}$. Their Coxeter diagrams are
shown in Table~\ref{t5}.

\end{lemma}

\begin{proof}
In fact, the lemma has been proved above. Here we show how to
verify the previous considerations without any geometry and
without referring to classifications from~\cite{Kap}
and~\cite{trun}. Since the procedure is very similar to the proof
of Lemma~\ref{ldim5-10}, we provide only a plan of necessary
computations without details.

Let $P$ be a compact hyperbolic Coxeter polytope $P$ with
Gale diagram $\D_{9}$. This Gale diagram contains an arc
$J=\lf 1,4,1\rf_2$, so the Coxeter diagram $\Sigma$ of $P$
consists of one of the diagrams $\Sigma_J$ presented in Item $2$
of Table~\ref{dugi-t} and two nodes $u_7$ and $u_8$ joined by a
dotted edge.

Choose one of three diagrams $\Sigma_J$.
Consider the diagram $\Sigma'=\l\Sigma_J, u_{7}\r$. It is connected,
contains a unique dotted edge, no multi-multiple edges, and its
determinant is equal to zero. So, we are able to find the weight
of the dotted edge joining $u_7$ with $\Sigma_J$ depending on
multiplicities of the remaining edges incident to $u_7$. The weight of
this edge should be greater than one. Of course, we must restrict
ourselves to the cases when non-dotted edges incident to $u_7$ do
not make any new Lann\'er diagram together with $\Sigma_J$. The
number of such cases is really small.

Further, we do the same for the diagram $\Sigma''=\l\Sigma_J, u_{8}\r$,
and we find all possible such diagrams together with the weight of
the dotted edge joining $u_8$ with $\Sigma_J$. Then we are left to
determine the weight of the dotted edge $u_7u_8$ for any pair of
diagrams $\Sigma'$ and $\Sigma''$. It occurs that this weight is
always greater than one.

Doing the procedure described above for all the three possible
diagrams $\Sigma_J$, we obtain the complete list of compact
hyperbolic Coxeter $5$-polytopes with $8$ facets with Gale diagram
$\D_{9}$. The computations completely confirm the result of
considerations previous to the lemma.

\end{proof}

In the remaining part of this section we show that Gale diagrams
$\D_5$-- $\D_8$ do not give rise to any Coxeter polytope.

\begin{lemma}
\label{ldim5-8}
There are no compact hyperbolic Coxeter polytopes
with Gale diagram $\D_{8}$.

\end{lemma}

\begin{proof}
Suppose that there exists a compact hyperbolic Coxeter polytope $P$ with
Gale diagram $\D_{8}$. This Gale diagram contains an arc $J=\lf 2,3,1\rf_2$.
According to Lemma~\ref{dugi} (Item $9$ of Table~\ref{dugi-t}) and Lemma~\ref{L},
the Coxeter diagram $\Sigma$ of $P$ consists of one of the nine
subdiagrams $\Sigma_J$ shown in Table~\ref{231l},
\begin{table}[!h]
\begin{center}
\caption{All possible diagrams $\Sigma_J$ for $J=\lf 2,3,1\rf_2$}
\label{231l}
\medskip
\psfrag{1}{\small $u_1$}
\psfrag{2}{\small $u_2$}
\psfrag{3}{\small $u_3$}
\psfrag{4}{\small $u_4$}
\psfrag{5}{\small $u_5$}
\psfrag{6}{\small $u_6$}
\epsfig{file=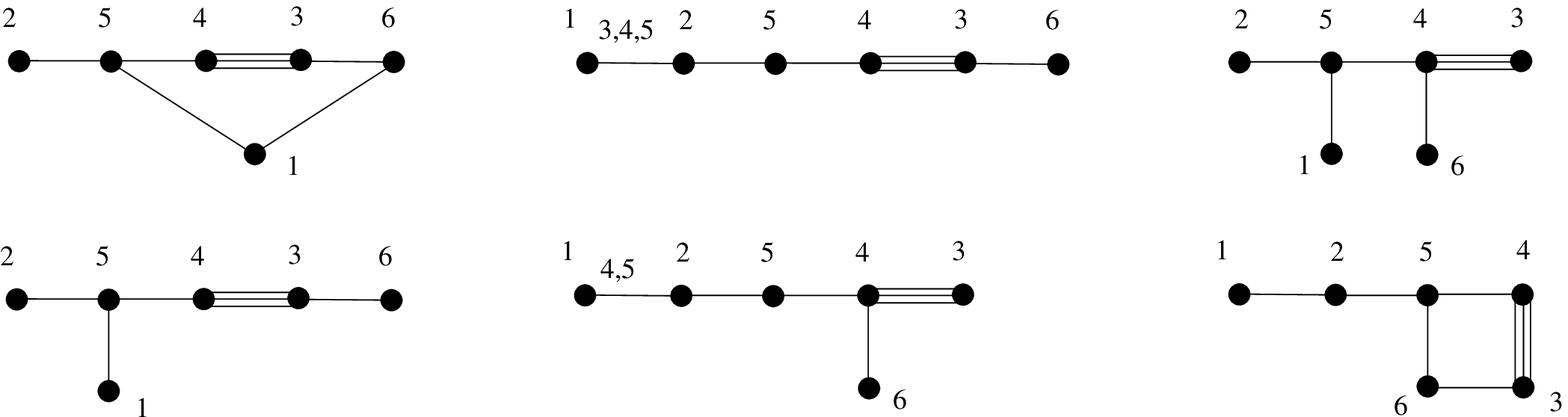,width=0.95\linewidth}
\end{center}
\end{table}
and two nodes $u_7$, $u_{8}$ joined by a dotted edge.
By Lemma~\ref{face}, the subdiagrams $\l u_{7},u_1,u_2\r$
and $\l u_6,u_8\r$ are Lann\'er diagrams, and no other Lann\'er
subdiagram of $\Sigma$ contains $u_7$ or $u_8$.

Consider the diagram $\Sigma'=\l\Sigma_J, u_{7}\r$. It is connected,
does not contain dotted edges, and its determinant is equal to
zero. Observe that the diagram $\l u_2,u_3,u_4,u_5\r$ is of the
type $H_4$. Since the diagram $\l u_7,u_2,u_3,u_4,u_5\r$ is
elliptic, this implies that $u_7$ is not joined with $\l
u_2,u_3,u_4,u_5\r$. Furthermore, notice that the diagram $\l
u_3,u_4,u_6\r$ is of the type $H_3$. Since the diagram $\l
u_7,u_3,u_4,u_6\r$ is elliptic, we obtain that $[u_7,u_6]=2$ or
$3$. Thus, for each of $9$ diagrams $\Sigma_J$ we have $2$
possibilities of attaching $u_7$ to $\Sigma_J\setminus u_1$.
Solving the equation $\det\Sigma'=0$, we compute the weight of
the edge $u_7u_1$. In all $18$ cases the result is not of the form
$\cos\frac{\pi}m$ for positive integer $m$, which proves the lemma.

\end{proof}

\begin{lemma}
\label{ldim5-7}
There are no compact hyperbolic Coxeter polytopes
with Gale diagram $\D_{7}$.

\end{lemma}

\begin{proof}
Suppose that there exists a hyperbolic Coxeter polytope $P$ with
Gale diagram $\D_{7}$. This Gale diagram contains an arc $J=\lf 1,3,1\rf_2$.
Therefore, the Coxeter diagram $\Sigma$ of $P$ contains one of the five
subdiagrams $\Sigma_J$, shown in Item $8$ of Table~\ref{dugi-t}.

On the other hand, $\Sigma$ contains a Lann\'er diagram $L$ of order
$3$ intersecting $\Sigma_J$. Denote by $u_1$ the intersection node
of $L$ and $\Sigma_J$, and denote by $u_6$ and $u_7$ the remaining
nodes of $L$. Since $L$ is connected, we may assume that $u_6$
attaches to $u_1$. Denote by $u_2$ the node of $\Sigma_J$
different from $u_1$ and contained in only one Lann\'er diagram of
order $4$, and denote by $u_3,u_4,u_5$ the nodes of $\Sigma_J$
contained in two Lann\'er diagrams of order $4$.

Consider the diagram $\Sigma_0=\l\Sigma_J, u_{6}\r\setminus u_2$.
It is connected, has order $5$, and contains a unique Lann\'er
diagram which is of order $4$. All such diagrams are listed
in~\cite[Lemma 3.8]{Ess2} (see the first two rows of Tabelle $3$,
the case $|{\cal N}_F|=1$, $|{\cal L}_F|=4$). We have reproduced
this list in Table~\ref{4pl1}.

Consider the diagram $\Sigma_1=\l\Sigma_J, u_{6}\r=\l\Sigma_J,
\Sigma_0\r$. Comparing the lists of possibilities for $\Sigma_J$
and $\Sigma_0$, it is easy to see that $\Sigma_1$ coincides with
one of the four diagrams listed in Table~\ref{1131l} (up to
permutation of indices $3,4$ and $5$).
\begin{table}[!h]
\begin{center}
\caption{All possibilities for diagram $\Sigma_1$, see Lemma~\ref{ldim5-7}}
\label{1131l}
\medskip
\psfrag{1}{\small $u_1$}
\psfrag{2}{\small $u_2$}
\psfrag{3}{\small $u_3$}
\psfrag{4}{\small $u_4$}
\psfrag{5}{\small $u_5$}
\psfrag{6}{\small $u_6$}
\psfrag{4,5}{\scriptsize $4,5$}
\epsfig{file=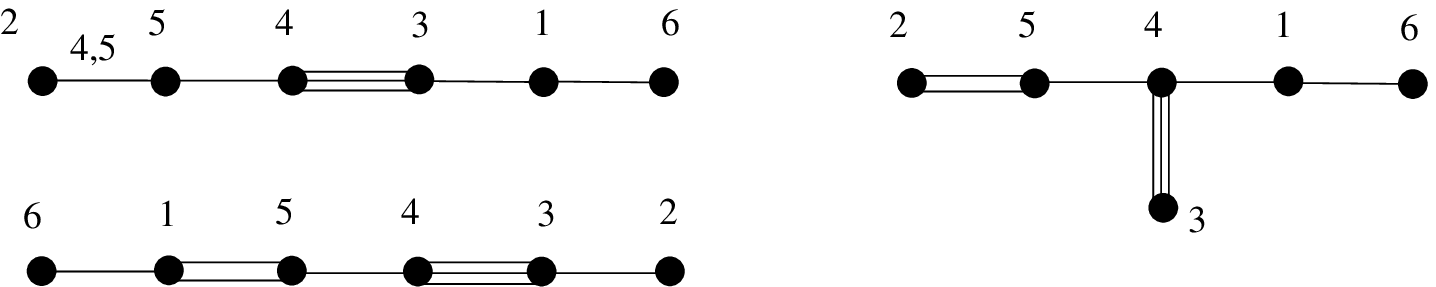,width=0.6\linewidth}
\end{center}
\end{table}
Now consider the diagram $\Sigma'=\l\Sigma_J, L\r=\l\Sigma_1,
u_7\r$. It is connected, does not contain dotted edges, its
determinant is equal to zero, and the only multi-multiple edge may
join $u_7$ and $u_6$. To reduce further computations notice, that
the diagram $\l u_7,u_3,u_4,u_5\r$ is elliptic, so $u_7$ does not
attach to $\l u_3,u_4\r$, and may attach to $u_5$ by simple edge
only. Moreover, since the diagrams $\l u_7,u_2,u_4,u_5\r$ and $\l
u_7,u_1,u_4,u_5\r$ are elliptic, $u_7$ is not joined with $u_5$.
Furthermore, since the diagrams $\l u_7,u_1,u_4,u_5\r$ and $\l
u_7,u_1,u_3,u_4\r$ are elliptic, $[u_7,u_1]=2$ or $3$. Considering
elliptic diagrams $\l u_7,u_2,u_4,u_5\r$ and $\l
u_7,u_2,u_3,u_4\r$, we obtain that $[u_7,u_2]$ is also at most
$3$. Then for all $4$ diagrams $\Sigma_1$ and all admissible
multiplicities of edges $u_7u_1$ and $u_7u_2$ we compute the
weight of the edge $u_7u_6$. We obtain exactly two diagrams
$\Sigma'$ where this weight is equal to $\cos\frac{\pi}m$ for some
positive integer $m$, these diagrams are shown in Fig.~\ref{2131l}.
\begin{figure}[!h]
\begin{center}
\psfrag{1}{\small $u_1$}
\psfrag{2}{\small $u_2$}
\psfrag{3}{\small $u_3$}
\psfrag{4}{\small $u_4$}
\psfrag{5}{\small $u_5$}
\psfrag{6}{\small $u_6$}
\psfrag{7}{\small $u_7$}
\psfrag{10}{\scriptsize $10$}
\epsfig{file=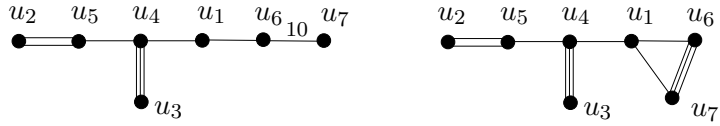,width=0.6\linewidth}
\caption{All possibilities for diagram $\Sigma'$, see Lemma~\ref{ldim5-7}}
\label{2131l}
\end{center}
\end{figure}
We are left to attach the node $u_8$ to $\Sigma'$.
Consider the diagram $\Sigma''=\Sigma\setminus u_2$. As usual,
it is connected, does not contain dotted edges, its determinant is
equal to zero, and the only multi-multiple edge that may appear is
$u_8u_7$. Furthermore, the diagram $\l u_3,u_4,u_1,u_6\r$ is of
the type $H_4$, and the diagram $\l u_8,u_3,u_4,u_1,u_6\r$ is
elliptic. Thus, $u_8$ does not attach to $\l u_3,u_4,u_1,u_6\r$.
The diagram $\l u_3,u_4,u_5\r$ is of the type $H_3$, and since
the diagram $\l u_8,u_3,u_4,u_5\r$ should be elliptic, this
implies that $[u_8,u_5]=2$ or $3$. Now for both diagrams
$\Sigma'\setminus u_2\subset\Sigma''$ we compute the weight of the
edge $u_8u_5$. In all four cases this weight is not equal to
$\cos\frac{\pi}m$ for any positive integer $m$, that finishes the
proof.

\end{proof}

\begin{lemma}
\label{ldim5-6}
There are no compact hyperbolic Coxeter polytope
with Gale diagram $\D_{6}$.

\end{lemma}

\begin{proof}
Suppose that there exists a hyperbolic Coxeter polytope $P$ with
Gale diagram $\D_{6}$.
The Coxeter diagram $\Sigma$ of $P$ consists of two Lann\'er diagrams
$L_1$ and $L_2$ of order $3$, and one Lann\'er diagram $L_3$ of
order $2$. Any two of these Lann\'er diagrams are joined in
$\Sigma$, and any subdiagram of $\Sigma$ not containing one of
these three diagrams is elliptic.

Consider the diagram $\Sigma_{12}=\l L_1,L_2\r$. Due
to~\cite[p. 239, Step 4]{Ess}, we have three cases:

(1) $L_1$ and $L_2$ are joined by two simple edges having a common
vertex, say in $L_2$;

(2) $L_1$ and $L_2$ are joined by a unique double edge;

(3) $L_1$ and $L_2$ are joined by a unique simple edge.\\

\noindent
We fix the following notation: $L_1=\l u_1,u_2,u_3\r$, $L_2=\l u_4,u_5,u_6\r$,
$L_3=\l u_7,u_8\r$, the only node of $L_2$ joined with $L_1$ is
$u_4$; $u_4$ is joined with $u_3$ and, in case (1), with $u_1$. We
may assume also that $u_7$ attaches to $L_1$, $u_4$ is joined to
$u_5$ in $L_2$, and $u_2$ is joined to $u_3$ in $L_1$.

\noindent
{\bf Case (1). }\ Since the diagrams $\l u_2,u_1,u_4\r$ and $\l
u_2,u_3,u_4\r$ are elliptic, $[u_2,u_1]$ and $[u_2,u_3]$ do not
exceed $5$. On the other hand, $\l u_1,u_2,u_3\r=L_1$ is a
Lann\'er diagram, so we may assume that $[u_2,u_1]=5$, and
$[u_2,u_3]=4$ or $5$. Now attach $u_7$ to $L_1$. If $u_7$ is
joined with $u_1$ or $u_2$, then the diagram $\l u_2,u_1,u_4\r$ is
not elliptic, and if $u_7$ is joined with $u_3$, then the diagram
$\l u_2,u_3,u_4\r$ is not elliptic, which contradicts 
Lemma~\ref{face}.

\noindent
{\bf Case (2). }\ It is clear that $[u_2,u_3]=[u_4,u_5]=3$, and
$u_7$ cannot be attached to $u_3$. Thus, $u_7$ is joined with
$u_1$ or $u_2$, which implies that $[u_2,u_1]\le 5$. Therefore,
$[u_1,u_3]=3$. So, the diagrams $\l u_1,u_3,u_4,u_5\r$ and $\l
u_2,u_3,u_4,u_5\r$ are of the type $F_4$. Therefore, if $u_7$
attaches $u_1$, then the diagram $\l u_7,u_1,u_3,u_4,u_5\r$ is not
elliptic, and if $u_7$ is joined with $u_2$, then the diagram $\l
u_7,u_2,u_3,u_4,u_5\r$ is not elliptic.

\noindent
{\bf Case (3). }\ The signature of $\Sigma_{12}$ is either $(5,1)$
or $(4,1,1)$. Thus, $\det\Sigma_{12}\le 0$. By Prop.~\ref{loc_product},
$\det(L_1,u_3)\det(L_2,u_4)\le\frac{1}{4}$. We may assume that
$|\det(L_1,u_3)|\le |\det(L_2,u_4)|$, in particular,
$|\det(L_1,u_3)|\le\frac{1}{2}$. By~\cite[Table 2]{Ess}, there are
only $6$ possibilities for $\l L_1,u_4\r$, we list them in
Table~\ref{l-14}.
\begin{table}[!h]
\begin{center}
\caption{All possibilities for diagram $\l L_1,u_4\r$, see Case $(3)$ of Lemma~\ref{ldim5-6}}
\label{l-14}
\medskip
\psfrag{1}{\small $u_1$}
\psfrag{2}{\small $u_2$}
\psfrag{3}{\small $u_3$}
\psfrag{4}{\small $u_4$}
\psfrag{7}{\scriptsize $7$}
\epsfig{file=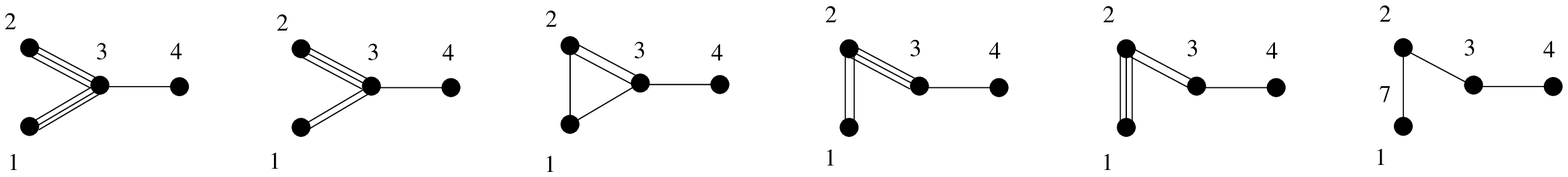,width=0.9\linewidth}
\end{center}
\end{table}

\noindent
For any of these six diagrams
$|\det(L_1,u_3)|\ge\frac{\sqrt{5}-1}8$. Thus,
$|\det(L_2,u_4)|\le\frac{1}{4}\frac{8}{\sqrt{5}-1}=\frac{2}{\sqrt{5}-1}$.
Notice that since the diagrams $\l u_3,u_4,u_5\r$ and $\l u_3,u_4,u_6\r$
are elliptic, $[u_4,u_5]$ and $[u_4,u_6]$ do not exceed $5$. Now,
since the local determinant is an increasing function of
multiplicities of the edges, it is not difficult to list all
Lann\'er diagrams $L_2=\l u_4,u_5,u_6\r$, such that
$[u_4,u_5],[u_4,u_6]\le 5$, and
$|\det(L_2,u_4)|\le\frac{2}{\sqrt{5}-1}$.
This list contains $17$ diagrams only.

Then, from $6\cdot 17=102$ pairs $(L_1,L_2)$ we list all pairs with
$\det(L_1,u_3)\det(L_2,u_4)\le\frac{1}{4}$. Each of these pairs
corresponds to a diagram $\Sigma_{12}$. After that, we attach to all
diagrams $\Sigma_{12}$ a node $u_7$ in the following way: $u_7$ is
joined with $L_1$ (and may be joined with $L_2$, too), and it does
not produce any new Lann\'er or parabolic diagram. It occurs that
none of obtained diagrams $\l \Sigma_{12},u_7\r$ has zero
determinant.

\end{proof}

\begin{lemma}
\label{ldim5-5}
There are no compact hyperbolic Coxeter polytopes
with Gale diagram $\D_{5}$.

\end{lemma}

\begin{proof}
Suppose that there exists a hyperbolic Coxeter polytope $P$ with
Gale diagram $\D_{5}$. The Coxeter diagram $\Sigma$ of $P$ consists
of one Lann\'er diagram $L_1$ of order $4$, and two
Lann\'er diagrams $L_2$ and $L_3$ of order $2$. Any two of these Lann\'er
diagrams are joined in $\Sigma$, and any subdiagram of $\Sigma$
not containing one of these three diagrams is elliptic.

We fix the following notation: $L_1=\l u_1,u_2,u_3,u_4\r$, $L_2=\l u_5,u_7\r$,
$L_3=\l u_6,u_8\r$, $u_5$ and $u_6$ attach to $L_1$.

Consider the diagram $\Sigma_0=\l L_1,u_5,u_{6}\r$.
It is connected, has order $6$, and contains a unique Lann\'er
diagram which is of order $4$. All such diagrams are listed
in~\cite[Lemma 3.8]{Ess2} (see Tabelle $3$, the first two rows of
page $27$, the case $|{\cal N}_F|=2$, $|{\cal L}_F|=4$). We have
reproduced this list in Table~\ref{4pl2}. The list contains five
diagrams, but we are interested in four of them: in the fifth one
only one of two additional nodes attaches to the
Lann\'er diagram. We list these four possibilities for $\Sigma_0$
in Table~\ref{4pl2u}.
\begin{table}[!h]
\begin{center}
\caption{All possibilities for diagram $\Sigma_0$, see Lemma~\ref{ldim5-5}}
\label{4pl2u}
\medskip
\psfrag{1}{\small $u_1$}
\psfrag{2}{\small $u_2$}
\psfrag{3}{\small $u_3$}
\psfrag{4}{\small $u_4$}
\psfrag{5}{\small $u_5$}
\psfrag{6}{\small $u_6$}
\epsfig{file=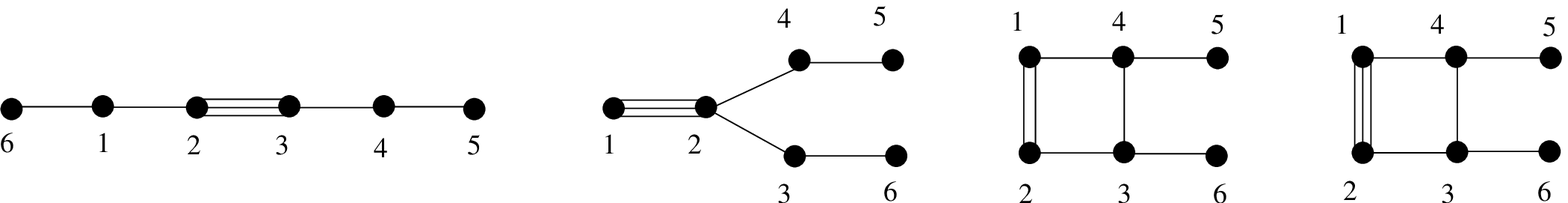,width=0.9\linewidth}
\end{center}
\end{table}

Now consider the diagram $\Sigma'=\l\Sigma_0,u_{7}\r$. It contains
a unique dotted edge $u_5u_7$. Since the diagram $\l
u_7,u_1,u_2,u_3,u_6\r$ is elliptic and the diagram $\l
u_1,u_2,u_3,u_6\r$ is of the type $H_4$ or $B_4$, $u_7$ is not
joined with $\l u_1,u_2,u_3\r$, and it may attach to $u_6$ if
$[u_1,u_2]=4$ only. It is easy to see that $[u_7,u_4]=2$ or $3$ in
all four cases. We obtain $9$ possibilities for attaching $u_7$ to
$\Sigma_0\setminus u_{5}$. For each of them we compute the weight
of the edge $u_5u_7$.

By symmetry, we may list all $9$ possibilities for the diagram
$\Sigma''=\l\Sigma_0,u_{8}\r$. Now we are left to compute the
weight of the edge $u_7u_8$ in $\Sigma$. Diagrams $\Sigma_0$ with
$[u_1,u_2]=5$ produce three possible diagrams $\Sigma$ each, and
the diagram $\Sigma_0$ with $[u_1,u_2]=4$ produces six possible
diagrams $\Sigma$ (we respect symmetry). In all these $15$ cases
the weight of the edge $u_7u_8$ is not of the form
$\cos\frac{\pi}m$ for positive integer $m$.

\end{proof}

\subsection{Dimension 4}

In dimension $4$ we must consider four Gale diagrams, namely $\D_{1}$-- $\D_{4}$.
Three of them, i.e. $\D_{1}$, $\D_{2}$ and $\D_{4}$, give rise to
Coxeter polytopes.

\begin{lemma}
\label{ldim4-1}
There are exactly three compact hyperbolic Coxeter polytopes
with Gale diagram $\D_{1}$. Their Coxeter diagrams are shown in
the third row of the second part of Table~\ref{t4}.

\end{lemma}

\begin{proof}
Let $P$ be a compact hyperbolic Coxeter polytope with
Gale diagram $\D_{1}$. The Coxeter diagram $\Sigma$ of $P$ consists
of one Lann\'er diagram $L_1$ of order $3$, and two
Lann\'er diagrams $L_2$ and $L_3$ of order $2$. Any two of these Lann\'er
diagrams are joined in $\Sigma$, and any subdiagram of $\Sigma$
containing none of these three diagrams is elliptic.

On the first sight, the considerations may repeat ones from the
proof of Lemma~\ref{ldim5-5}. However, there is a small
difference: the number of Lann\'er diagrams of order $3$ is
infinite. Thus, at first we must bound the multiplicities of the
edges of the Lann\'er diagram of order $3$.

We fix the following notation: $L_1=\l u_1,u_2,u_3\r$, $L_2=\l u_5,u_6\r$,
$L_3=\l u_4,u_7\r$, $u_4$ and $u_5$ attach to $L_1$. We may also
assume that $u_4$ attaches to $u_3$.

Since the diagrams $\l u_1,u_3,u_4\r$ and $\l u_2,u_3,u_4\r$
should be elliptic, the edges $u_3u_1$ and $u_3u_2$ are not
multi-multiple. We consider two cases: $u_1$ or $u_2$ is either
joined with $\l u_4,u_5,u_6,u_7\r$ or not.\\

\noindent
{\bf Case 1: }$u_1$ and $u_2$ are not joined with $\l
u_4,u_5,u_6,u_7\r$. In particular, this is true if the edge
$u_1u_2$ is multi-multiple. Then $u_5$ attaches to $u_3$.
Since the diagrams $\l u_1,u_3,u_4,u_5\r$ and $\l
u_2,u_3,u_4,u_5\r$ are elliptic, $[u_3,u_1]$ and $[u_3,u_2]$ do
not exceed $3$, $[u_3,u_4]=[u_3,u_5]=3$, and $[u_4,u_5]=2$. We may
assume that $[u_3,u_1]=3$, and $[u_3,u_2]=2$ or $3$.

Consider the diagram $\Sigma'=\l L_1, L_2,u_4\r=\Sigma\setminus
u_7$. We know that $u_6$ is joined with $u_5$ by a dotted edge,
and $u_6$ does not attach to $u_1$ and $u_2$. Furthermore, since
the diagram $\l u_1,u_3,u_4,u_6\r$ is elliptic, $[u_6,u_3]\le 3$
and $[u_6,u_4]\le 4$. By the same reason, either $[u_6,u_3]$ or
$[u_6,u_4]$ is equal to $2$. Thus, we have four possibilities to
attach $u_6$ to $u_3$ and $u_4$.

Denote by $d$ the weight of the dotted edge $u_5u_6$, and compute
the local determinant $\det\left(\l u_3,u_4,u_5,u_6\r,u_3\right)$
for all four diagrams $\l u_3,u_4,u_5,u_6\r$ as a function of $d$.

\smallskip

\noindent
{\bf Case 1.1: } $[u_6,u_4]\ne 2$. In this case $\det\left(\l
u_3,u_4,u_5,u_6\r,u_3\right)$ equals either
$\frac{12d^2+4d-5}{4(4d^2-3)}$ (when $[u_6,u_4]=3$) or
$\frac{6d^2+2\sqrt{2}d-1}{4(2d^2-1)}$ (when $[u_6,u_4]=4$).
Both expressions decrease in the ray $[1,\infty)$, so the maximal
values are $11/4$ and $(5+2\sqrt{2})/4$ respectively.
Now recall that $\det\Sigma'=0$, so by Prop.~\ref{loc_sum} we
have $\det(L_1,u_3)=1-\det\left(\l u_3,u_4,u_5,u_6\r,u_3\right)$.
Therefore, $|\det(L_1,u_3)|$ is bounded from above by $7/4$ or
$(1+2\sqrt{2})/4$ if $[u_6,u_4]=3$ or $[u_6,u_4]=4$ respectively.
Since $|\det(L_1,u_3)|$ is an increasing function on $[u_1,u_2]$,
an easy check shows that $[u_1,u_2]$ is bounded by $10$ or $8$
respectively. So, in both cases we have finitely many
possibilities for $L_1$.

Further considerations follow ones from
Lemma~\ref{ldim5-5}. We list all possible $\Sigma'$ together with
the weight of the dotted edge $u_5u_6$ (which may be computed from
the equation $\det\Sigma'=0$), then we list all possible
diagrams $\Sigma''=\l L_1, L_3,u_5\r=\Sigma\setminus u_6$ in a
similar way.
After that for all pairs $\left(\Sigma',\Sigma''\right)$ (with the
same $L_1$) we compute the weight of the edge $u_6u_7$. It occurs
that in all cases the weight is not of the form $\cos\frac{\pi}m$
for positive integer $m$.

\smallskip

\noindent
{\bf Case 1.2: } $[u_6,u_4]=2$. In this case $\det\left(\l
u_3,u_4,u_5,u_6\r,u_3\right)$ equals either
$\frac{3d^2-2}{4(d^2-1)}$ (when $[u_6,u_3]=2$) or
$\frac{3d-1}{4(d-1)}$ (when $[u_6,u_3]=3$). These tend to $\infty$
when $d$ tends to $1$, so we do not obtain any bound for
$[u_1,u_2]$.

Let $m_{12}=[u_1,u_2]$, $m_{23}=[u_2,u_3]$, and let
$m_{36}=[u_3,u_6]$. Notice that $m_{23},m_{36}=2$ or $3$.
Define also $c_{12}=\cos(\pi/m_{12})$.
We compute the weight of the edge $u_5u_6$ as a function
$d(m_{12},m_{23},m_{36})$ of $m_{12}$, $m_{23}$ and $m_{36}$.
Solving the equation $\det\Sigma'=0$, we see that
\begin{multline*}
\qquad d(m_{12},2,2)=\sqrt{\frac{2c^2_{12}-1}{2c^2_{12}-2}};\qquad\qquad
d(m_{12},3,2)=\sqrt{\frac{2c_{12}}{3c_{12}-1}};\\
d(m_{12},2,3)={\frac{c^2_{12}}{3c^2_{12}-2}};\qquad\qquad
d(m_{12},3,3)={\frac{c_{12}+1}{3c_{12}-1}}.\qquad
\end{multline*}
Consider the diagram $\Sigma$. According to Case~1.1, we may assume
that $[u_5,u_7]=2$. Since $L_2$ and $L_3$ are joined in
$\Sigma$, $[u_6,u_7]\ne 2$. On the other hand, the diagram $\l
u_3,u_6,u_7\r$ is elliptic. Thus, either $[u_3,u_6]$ or
$[u_3,u_7]$ equals $2$. By symmetry, we may assume that
$[u_3,u_7]=2$. We also know how the weight of the edge $u_4u_7$
depends on  $m_{12}$ and $m_{23}$.

Now we are able to compute the weight of the dotted edge $u_4u_7$ as a
function $w(m_{12},m_{23},m_{36})$ of $m_{12}, m_{23}$ and
$m_{36}$. For that we simply solve the equation $\det\Sigma=0$.
Notice that since $L_1$ is a Lann\'er diagram, $m_{12}\ge 7$ when
$m_{23}=2$, and $m_{12}\ge 4$ when $m_{23}=3$.
We obtain:

$\bullet$
$w(m_{12},2,2)={\frac{\textstyle 1-c^2_{12}}{\textstyle 3c^2_{12}-2}}$ is a decreasing
function of $m_{12}$ while $m_{12}\ge 7$, and $w(7,2,2)<1/2$;

$\bullet$
$w(m_{12},2,3)={\frac{\textstyle 2(1-c^2_{12})\sqrt{2c^2_{12}}}%
{\textstyle (3c^2_{12}-2)^{3/2}}}$
is a decreasing function of $m_{12}$ while $m_{12}\ge 7$,
$w(9,2,3)<1/2$, and $w(m_{12},2,3)\ne \cos(\pi/m)$ when $m_{12}=7$
or $8$;

\medskip
$\bullet$
$w(m_{12},3,2)={\frac{\textstyle 1-c_{12}}{\textstyle 3c_{12}-1}}$ is a decreasing
function of $m_{12}$ while $m_{12}\ge 4$, and $w(4,2,2)<1/2$;

$\bullet$
$w(m_{12},3,3)={\frac{\textstyle 2(1-c_{12})\sqrt{2c_{12}}}{\textstyle (3c_{12}-1)^{3/2}}}$
is a decreasing function of $m_{12}$ while $m_{12}\ge 4$,
$w(5,3,3)<1/2$, and $w(4,3,3)\ne \cos(\pi/m)$.

\medskip
This finishes considerations of Case~1.\\

\noindent
{\bf Case 2: }either $u_1$ or $u_2$ is joined with $\l
u_4,u_5,u_6,u_7\r$. In particular, this implies that $L_1$
contains no multi-multiple edges, so we deal with a finite number
of possibilities for $L_1$ only. This list contains $11$
Lann\'er diagrams of order $3$.
Using that list, it is not too difficult to list all the diagrams
$\Sigma_0=\l L_1,u_4,u_{5}\r$. This list contains $19$ diagrams,
we present them in Table~\ref{131u}.
\begin{table}[!h]
\begin{center}
\caption{All possibilities for diagram $\Sigma_0$, see Case $2$ of Lemma~\ref{ldim4-1}}
\label{131u}
\medskip
\psfrag{1}{\small $u_1$}
\psfrag{2}{\small $u_2$}
\psfrag{3}{\small $u_3$}
\psfrag{4}{\small $u_4$}
\psfrag{5}{\small $u_5$}
\psfrag{6}{\small $u_6$}
\psfrag{7}{\small $u_7$}
\epsfig{file=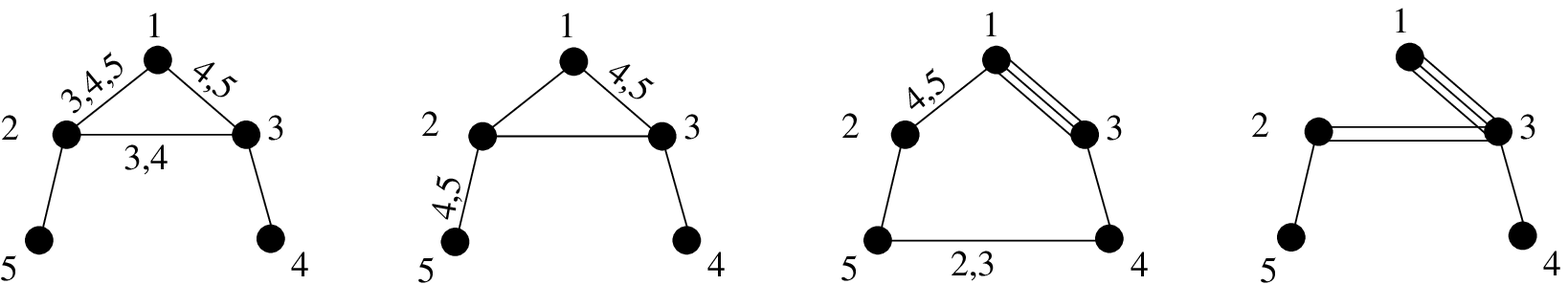,width=0.9\linewidth}
\end{center}
\end{table}
Now we follow the proof of Lemma~\ref{ldim5-5}. Choose one of $19$
diagrams $\Sigma_0$, and consider the diagram
$\Sigma'=\l\Sigma_0,u_{6}\r$. It contains a unique dotted edge
$u_5u_6$, and that is the only Lann\'er diagram in $\Sigma'$ containing $u_6$.
We have a finite number of possibilities to attach $u_6$ to
$\Sigma_0\setminus u_5$. For each of them we compute the weight of
the edge $u_5u_6$.

Similarly, we list all possibilities for the diagram
$\Sigma''=\l\Sigma_0,u_{7}\r$. Now we are left to compute the
weight of the edge $u_6u_7$ in $\Sigma$. A computation shows that
the weight is of the form $\cos\frac{\pi}m$ only for the diagrams
listed in Table~\ref{322u-r}.
\begin{table}[!h]
\begin{center}
\caption{Coxeter diagrams of Coxeter polytopes with Gale diagram $\D_{1}$}
\label{322u-r}
\medskip
\psfrag{1}{\small $u_1$}
\psfrag{2}{\small $u_2$}
\psfrag{3}{\small $u_3$}
\psfrag{4}{\small $u_4$}
\psfrag{5}{\small $u_5$}
\psfrag{6}{\small $u_6$}
\psfrag{7}{\small $u_7$}
\psfrag{2,3}{\small $2,3$}
\psfrag{4,5}{\small $4,5$}
\psfrag{3,4,5}{\small $3,4,5$}
\epsfig{file=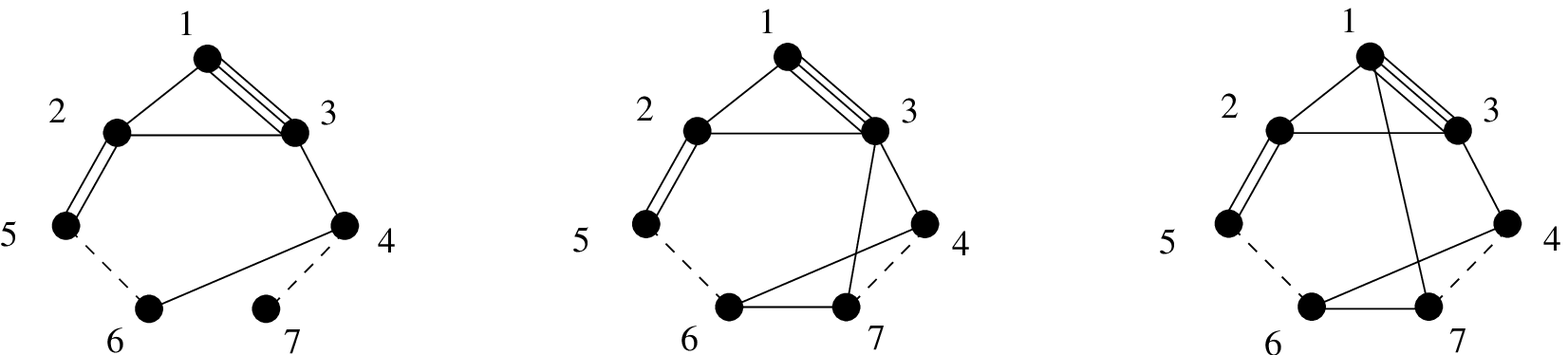,width=0.9\linewidth}
\end{center}
\end{table}
To verify that these diagrams correspond to polytopes, we need to
assign weights to the dotted edges. We assign a weight
$\sqrt{2}\frac{\sqrt{5}+1}{4}$ to all edges $u_5u_6$, and weights
$\frac{\sqrt{15(5+\sqrt{5})}}{10}$, $\frac{5+3\sqrt{5}}{10}$ and
$\frac{3+\sqrt{5}}{4}$ to the edge $u_4u_7$ on the left, middle
and right diagrams respectively. A direct calculation shows that
the diagrams have signature $(4,1,2)$.

\end{proof}

\begin{lemma}
\label{ldim4-2}
There are $29$ compact hyperbolic Coxeter polytopes with
Gale diagram $\D_{2}$. Their Coxeter diagrams are
shown in the first part of Table~\ref{t4} and in the first three
rows of the second part of the same table.

\end{lemma}

\begin{proof}
The proof is identical to one which concerns the diagram
$\D_9$ (see Lemma~\ref{ldim5-9}). The combinatorial type of
polytope defined by Gale diagram $\D_{2}$ is twice truncated
$4$-simplex. Any such Coxeter polytope may be obtained by gluing
one or two prisms to a twice truncated $4$-simplex with
orthogonality conditions described before Lemma~\ref{ldim5-9}.
Such simplices were classified by Schlettwein in~\cite{trun}, they
appear as right ones in rows $1$, $2$, and $4$ of the first part
of Table~\ref{t4}, and in rows $1$ and $2$ of the second part.
The prisms were classified by Kaplinskaja in~\cite{Kap}.

For each twice truncated simplex from the list of Schlettwein
we find all the prisms that have "right" base congruent to one of
"right" facets of the truncated simplex, and glue these prisms to
the truncated simplex. The result is presented in Table~\ref{t4}.

The verification of the result above by computations is
completely identical to the proof of Lemma~\ref{ldim5-9}.
We only need to replace an arc $J=\lf 1,4,1\rf_2$ from $\D_9$ by an arc
$J=\lf 1,3,1\rf_2$, and refer to Item $8$ of
Table~\ref{dugi-t} instead of Item $2$.

\end{proof}

\begin{lemma}
\label{ldim4-3}
There are no compact hyperbolic Coxeter polytopes
with Gale diagram $\D_{3}$.

\end{lemma}

\begin{proof}
Suppose that there exists a hyperbolic Coxeter polytope $P$ with
Gale diagram $\D_{3}$. The Coxeter diagram $\Sigma$ of $P$ consists
of one Lann\'er diagram $L_1=\l u_1,u_2,u_3,u_4\r$ of order $4$,
two Lann\'er diagrams $L_2=\l u_6,u_1,u_2\r$ and $L_3=\l
u_3,u_4,u_5\r$ of order $3$, and two Lann\'er diagrams $\l
u_6,u_7\r$ and $\l u_7,u_5\r$ of order $2$.

Consider the diagram $\Sigma'=\l L_1,L_2,L_3\r=\Sigma\setminus u_7$.
It is connected, has order $6$, and contains no dotted edges.
We may also assume that $u_5$ attaches to $u_4$.
Clearly, any multi-multiple edge that may appear in $\Sigma'$
belongs to $L_2$ or $L_3$ and does not belong to $L_1$. We
consider two cases: either $\Sigma'$ contains multi-multiple edges
or not.

Suppose that $\Sigma'$ contains no multi-multiple edges.
Then we have $9$ possibilities for $L_2$, and $9$ possibilities for
$L_3$. For each of $81$ pairs (or $45$ in view of symmetry) we
join nodes of $L_2$ with nodes of $L_3$ in all possible ways ($9$
edges, $4$ possibilities for each of them, from empty to triple
one). We are looking for diagrams satisfying the following
conditions: the determinant should vanish, there are no parabolic
subdiagrams, and the diagram contains a unique new Lann\'er
diagram, which has order $4$.
A computer check (which is not very short) shows that only $39$ obtained
diagrams have zero determinant, and only $11$ of them contain
Lann\'er diagrams of order $4$. However, each of them contains
some new Lann\'er diagram of order $3$. Therefore, none of them
may be considered as $\Sigma'$.

Now suppose that $\Sigma'$ contains at least one multi-multiple edge.
We may assume that $u_4u_5$ is multi-multiple.
In this case $u_4$ must be a {\it leaf} of $L_1$, i.e. it should have
valency one in $L_1$. Indeed, if $u_4$ is joined with two vertices
$v,w\in L_1$, then both diagrams $\l u_5,u_4,v\r$ and $\l
u_5,u_4,w\r$ are not elliptic, which is impossible. Thus, $L_1$ is
not a cycle, so we have $4$ possibilities for $L_1$ only (see
Table~\ref{lan}). In Table~\ref{tl4-22u} we list all possible diagrams $L_1$
together with all possible numerations of nodes. A numeration
should satisfy the following properties: $u_4$ is a leaf, and
$u_3$ is a unique neighbor of $u_4$. We consider numerations up to
interchange of $u_1$ and $u_2$.
\begin{table}[!h]
\begin{center}
\caption{Numberings of vertices of Lann\'er diagrams of order $4$ without cycles}
\label{tl4-22u}
\medskip
\psfrag{1}{\small $u_1$}
\psfrag{2}{\small $u_2$}
\psfrag{3}{\small $u_3$}
\psfrag{4}{\small $u_4$}
\psfrag{5}{\small $u_5$}
\psfrag{6}{\small $u_6$}
\psfrag{7}{\small $u_7$}
\psfrag{2a}{\small (2a)}
\psfrag{4a}{\small (4a)}
\psfrag{2b}{\small (2b)}
\psfrag{4b}{\small (4b)}
\psfrag{1a}{\small (1)}
\psfrag{3a}{\small (3)}
\epsfig{file=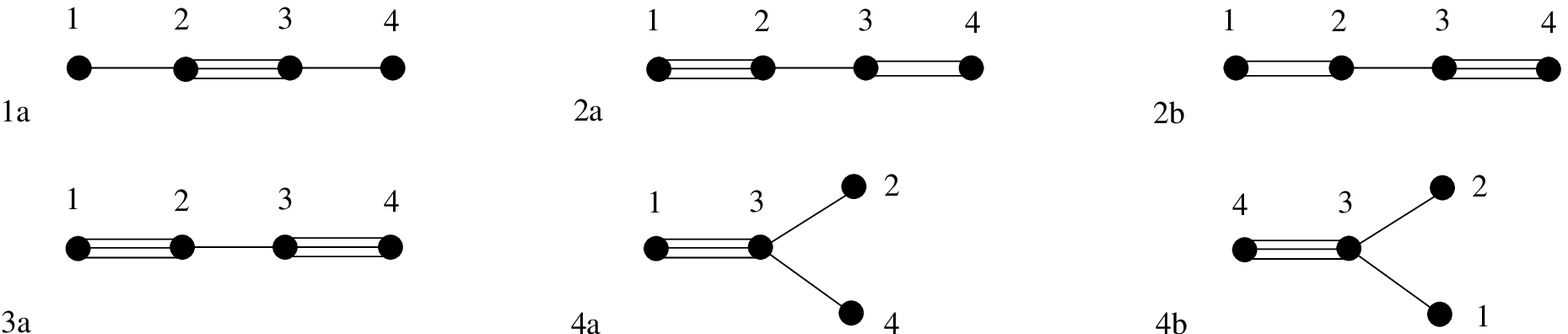,width=0.9\linewidth}
\end{center}
\end{table}

Consider $6$ diagrams case-by-case. For all of them we claim that
$u_5$ and $u_4$ do not attach to $L_2=\l u_1,u_2,u_6\r$: this is
because the edge $u_4u_5$ is multi-multiple.\\

\noindent
{\bf Diagram (1). }Since the diagram $\l u_1,u_2,u_3,u_5\r$ is
elliptic, $u_5$ is not joined with $u_3$. Furthermore, since the
diagram $\l u_6,u_2,u_3,u_4\r$ is elliptic, $u_6$ is not joined
with $\l u_2,u_3\r$. Therefore, $[u_6,u_2]=2$, so $[u_6,u_1]\ge
7$. Applying Prop.~\ref{loc_product}, we see that
$\det(L_2,u_2)\det(L_3,u_3)=\cos^2(\pi/5)$. An easy calculation
shows that the inequality $[u_6,u_1]\ge 7$ implies that
$[u_4,u_5]\le 10$. By symmetry, $[u_6,u_1]\le 10$, too. We are
left with a finite (and very small) number of possibilities for
$\Sigma'$. For none of them $\det\Sigma'=0$.\\

\noindent
{\bf Diagrams (2a), (2b) and (3). }Since the diagram $\l
u_1,u_2,u_3,u_5\r$ is elliptic, $[u_3,u_5]\le 3$. Since the
diagram $\l u_6,u_2,u_3,u_4\r$ is elliptic, $u_6$ is not joined
with $u_3$, and $[u_6,u_2]\le 3$, so $[u_6,u_1]\ge 3$. Applying
Prop.~\ref{loc_product}, we have $\det(L_2,u_2)\det(L_3,u_3)=1/4$.
By assumption, $[u_4,u_5]\ge 6$, which implies the inequality
$|\det(L_3,u_3)|\ge |\d(2,4,6)|=1$. Thus, $|\det(L_2,u_2)|\le
1/4$. But since $[u_1,u_2]\ge 4$ and $[u_6,u_2]\ge 3$, either
$|\det(L_2,u_2)|\ge |\d(2,4,5)|=1/\sqrt{5}>1/4$ or
$|\det(L_2,u_2)|\ge |\d(3,4,3)|=\sqrt{2}/3>1/4$, so we come to
a contradiction.\\

\noindent
{\bf Diagram (4a). }Since the diagram $\l u_6,u_1,u_3\r$ is
elliptic, $[u_6,u_1]\le 3$. On the other hand, $L_2=\l
u_1,u_2,u_6\r$ is a Lann\'er diagram, so $[u_6,u_2]\ge 7$. This
implies that $\l u_6,u_2,u_3\r$ is a Lann\'er diagram, which is
impossible.\\

\noindent
{\bf Diagram (4b). }Since the diagram $\l u_6,u_2,u_3,u_4\r$ is
elliptic, $[u_6,u_2]\le 3$. Hence, $[u_6,u_1]\ge 7$, and $\l
u_6,u_1,u_3\r$ is a Lann\'er diagram. This contradiction completes
the proof of the lemma.

\end{proof}

\begin{lemma}
\label{ldim4-4}
There are exactly eight compact hyperbolic Coxeter $4$-polytopes with $7$
facets with Gale diagram $\D_{4}$. Their Coxeter diagrams are
shown in the bottom of the second part of Table~\ref{t4}.

\end{lemma}

\begin{proof}
Let $P$ be a hyperbolic Coxeter polytope with  Gale diagram $\D_{4}$.
The Coxeter diagram $\Sigma$ of $P$ contains
two Lann\'er diagrams $L_1=\l u_1,u_2,u_3\r$ and $L_2=\l
u_3,u_4,u_5\r$ of order $3$, a dotted edge $u_6u_7$, and other two
Lann\'er diagrams $L_3=\l u_1,u_2,u_6\r$ and $L_4=\l
u_7,u_4,u_5\r$ of order $3$. Any subdiagram of $\Sigma$ 
containing none of these five diagrams is elliptic.
Since $L_3$ and $L_4$ are connected, we may assume that
$u_6$ attaches to $u_2$, and $u_7$ attaches to $u_5$.

Consider the diagram $\Sigma'=\l L_3,L_1,L_2\r=\Sigma\setminus u_7$.
Clearly, the only multi-multiple edges that may appear in $\Sigma'$ are
$u_1u_2$, $u_6u_2$, $u_6u_1$, and $u_4u_5$.

At first, suppose that the edge $u_6u_2$ is multi-multiple. Then
$\l u_6,u_2\r$ is not joined with $\l u_3,u_4,u_5\r=L_2$.
In particular, $[u_2,u_3]=2$, so $[u_1,u_3]\ne 2$. Thus,
$[u_6,u_1]$ is also equal to $2$.
Furthermore, since diagrams $\l u_1,u_3,u_4\r$ and $\l
u_1,u_3,u_5\r$ are elliptic, $[u_3,u_4], [u_3,u_5]$ and
$[u_1,u_3]\le 5$. Therefore, since $\l u_1,u_2,u_3\r=L_1$ is a
Lann\'er diagram, $[u_1,u_2]\ge 4$. Now suppose that $[u_1,u_4]\ne
2$. Then $[u_3,u_4]=2$, so $[u_4,u_5]\ge 4$, and the diagram $\l
u_2,u_1,u_4,u_5\r$ is not elliptic, which is impossible. The
contradiction shows that $[u_1,u_4]= 2$. Similarly, $[u_1,u_5]=
2$. Consequently, the diagram $\Sigma'$ looks like the diagram
shown in Fig.~\ref{1212k},
\begin{figure}[!h]
\begin{center}
\psfrag{1}{\small $u_1$}
\psfrag{2}{\small $u_2$}
\psfrag{3}{\small $u_3$}
\psfrag{4}{\small $u_4$}
\psfrag{5}{\small $u_5$}
\psfrag{6}{\small $u_6$}
\psfrag{>6}{\scriptsize $\ge \! 6$}
\psfrag{>4}{\scriptsize $\ge \!4$}
\psfrag{3,4,5}{\scriptsize $3,\!4,\!5$}
\psfrag{<5}{\scriptsize $\le \! 5$}
\psfrag{45}{\footnotesize $m_{45}$}
\epsfig{file=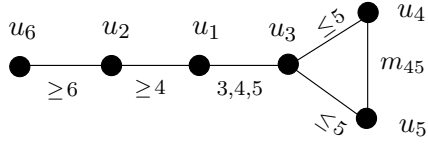,width=0.35\linewidth}
\caption{A diagram $\Sigma'$, see Lemma~\ref{ldim4-4}}
\label{1212k}
\end{center}
\end{figure}
where $m_{45}=[u_4,u_5]$. Now we may apply
Prop.~\ref{loc_product}:
$\det(L_3,u_1)\det(L_2,u_3)=\cos^2(\pi/m_{13})$, where
$m_{13}=[u_1,u_3]$.
Notice that since $[u_1,u_2]\ge 4$ and $[u_2,u_6]\ge 6$, we have
$|\det(L_3,u_1)|\ge|\d(2,4,6)|=1$.

If $m_{13}=4$ or $5$, we obtain that $[u_3,u_4],[u_3,u_5]\le 3$,
which implies $[u_4,u_5]= 7$ in view of
$|\det(L_2,u_3)|\le\cos^2(\pi/m_{13})$. Thus,
$|\det(L_2,u_3)|\ge|\d(2,3,7)|$. This implies that
$|\det(L_3,u_1)|\le\cos^2(\pi/5)/|\d(2,3,7)|$. An easy calculation
shows that in this case $[u_2,u_6]\le 7$, $[u_2,u_1]\le 6$. Then
we check the finite (small) number of possibilities for $\Sigma'$
and see that none of them has determinant equal to zero.

If $m_{13}=3$, then $[u_2,u_1]\le 7$. Therefore,
$|\det(L_3,u_1)|\ge|\d(2,6,7)|$. Hence,
$|\det(L_2,u_3)|\le\cos^2(\pi/3)/|\d(2,6,7)|$, but such $L_2$ does
not exist.

The contradiction shows that the edge $u_6u_2$ is not
multi-multiple. Similarly, the edges $u_6u_1$, $u_7u_5$, and
$u_7u_4$ of $\Sigma$ are not multi-multiple either.
Thus, the only edges that
may be multi-multiple in $\Sigma$ are $u_4u_5$ and $u_1u_2$.

Consider again the diagram $\Sigma'$ and suppose that the diagram $\l
u_4,u_5\r$ is not joined with $\l u_1,u_2,u_6\r$. In particular,
this holds if at least one of the edges $u_4u_5$ and $u_1u_2$ is
multi-multiple. We may apply
Prop.~\ref{loc_sum}:
$$\det(\l L_3,L_1\r,u_3)+\det(L_2,u_3)=1$$
By definition,
$$\det(\l L_3,L_1\r,u_3)=\det\l
L_3,L_1\r/\det(L_3)$$
We use a very rough bound: $|\det\l L_3,L_1\r|<16$
since it is a determinant of a $4\times 4$ matrix with entries
between $-1$ and $1$, and $|\det(L_3)|\ge |3/4-\cos^2(\pi/7)|=|\det(\L_{2,3,7})|$, 
since $\det(\L_{2,3,7})$ is maximal among all determinants of Lann\'er diagrams of 
order $3$. This bound implies
$$|\det(L_2,u_3)|\le 1+|\det(\l L_3,L_1\r,u_3)|\le 1+\frac{16}{|3/4-\cos^2(\pi/7)|}<261$$
Now an easy computation shows that $[u_4,u_5]\le 101$. Considering a diagram
$\Sigma''=\l L_1,L_2,L_4\r=\Sigma\setminus u_6$ in a similar way,
we obtain that $[u_1,u_2]\le 101$, too, and we are left with a
finite number of possibilities for $\Sigma'$ (and for $\Sigma''$).
We list all diagrams $L_2$ (less that $1000$ possibilities) and
all possible diagrams $\l L_3,L_1\r$ (less that $10000$
possibilities), and find all pairs such that $\det(\l
L_3,L_1\r,u_3)+\det(L_2,u_3)=1$, there are about $50$ such pairs.
Therefore, we obtain a complete list of possibilities for $\Sigma'$ (and for
$\Sigma''$). Then we look for unordered pairs
($\Sigma',\Sigma''$), such that the diagrams coincide on their
intersection, i.e. a subdiagram $\l
L_1,L_2\r\subset\Sigma'$ coincides with a subdiagram
$\l L_1,L_2\r\subset\Sigma''$. There are only $8$ such pairs, all
them give rise to Coxeter diagrams of Coxeter polytopes. The
diagrams are shown in the bottom of the second part of
Table~\ref{t4}. The weight of the dotted edge is equal to
$\sqrt{2}\cos(\pi/8)$ for the two last diagrams, is equal to
$(\sqrt{5}+1)/2$ for the three diagrams in the second row from the bottom, and is
equal to $1+\sqrt{2}$ for the three diagrams in the third row from the bottom.

Now suppose that the diagram $\l u_4,u_5\r$ is joined with
$\l u_1,u_2,u_6\r$. This implies that $\Sigma$ does not contain
multi-multiple edges, so we have a finite number of possibilities
for the diagrams $\Sigma'$ and $\Sigma''$. A computation shows
that we do not obtain any polytope in this way.

\end{proof}

The result of the considerations above is presented below.
Recall that there are no $7$-dimensional polytopes with $10$ facets.

\begin{table}[!h]
\caption{$8$-dimensional polytope with $11$ facets}
\label{t8}
\medskip
\begin{center}
\psfrag{d}{}
\epsfig{file=./pic/8/8rec.eps,width=0.51\linewidth}
\end{center}
\end{table}

\vspace{-20pt}

\begin{table}[!h]
\caption{$6$-dimensional polytopes with $9$ facets}
\label{t6}
\medskip
\begin{center}
\psfrag{10}{\footnotesize 10}
\psfrag{r1}{}
\psfrag{r}{}
\psfrag{d2}{}
\psfrag{d3}{}
\epsfig{file=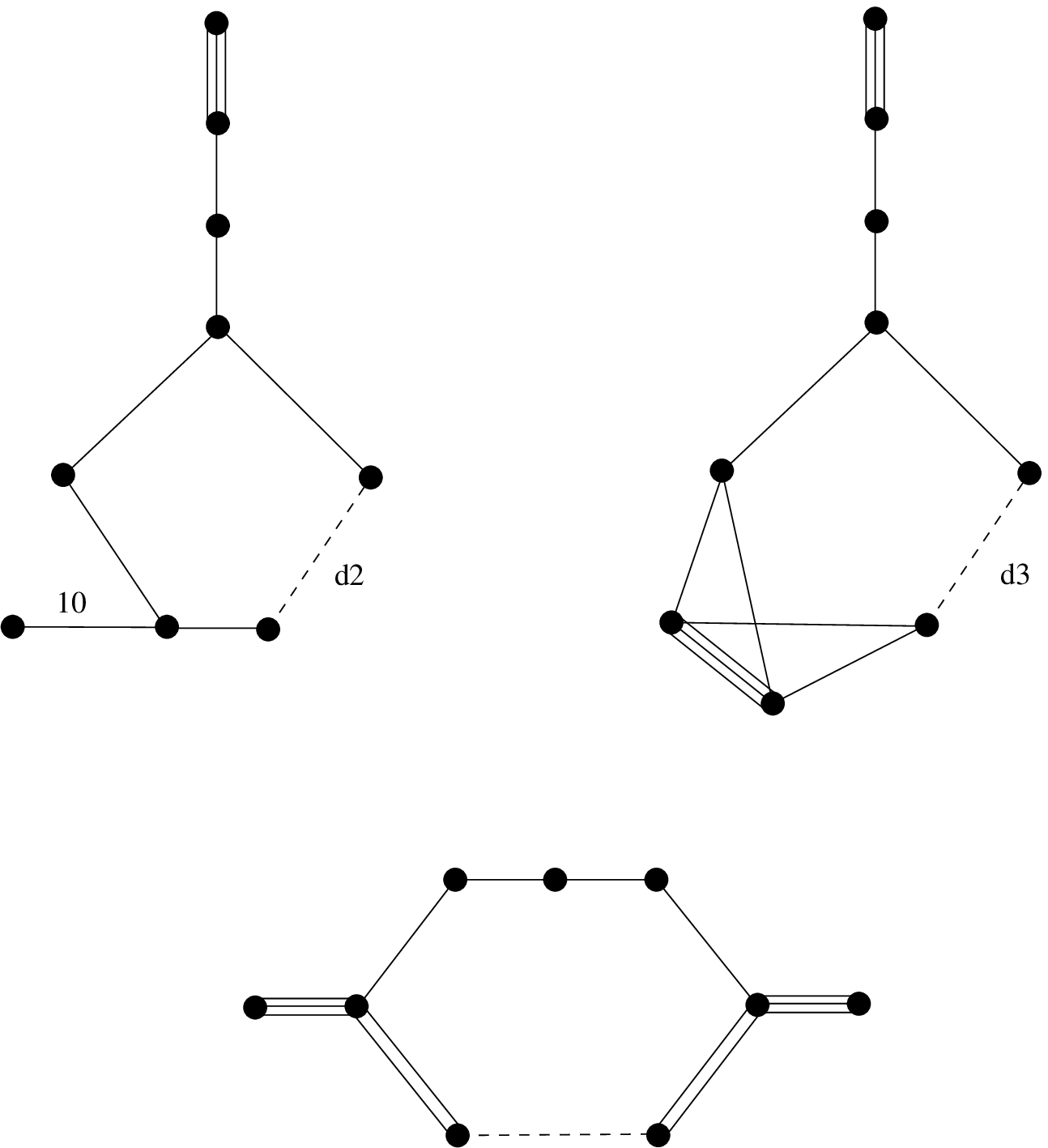,width=0.5\linewidth}
\end{center}
\end{table}

\pagebreak

\begin{table}[!h]
\caption{$5$-dimensional polytopes with $8$ facets}
\label{t5}
\medskip
\begin{center}
\psfrag{d11}{}
\psfrag{d12}{}
\psfrag{d13}{}
\psfrag{d21}{}
\psfrag{d22}{}
\psfrag{d23}{}
\psfrag{d31}{}
\psfrag{d32}{}
\psfrag{d33}{}
\psfrag{d41}{}
\psfrag{d42}{}
\psfrag{d43}{}
\psfrag{d51}{}
\psfrag{d52}{}
\psfrag{d53}{}
\psfrag{d61}{}
\psfrag{d62}{}
\psfrag{d63}{}
\psfrag{d71}{}
\psfrag{d72}{}
\psfrag{d73}{}
\psfrag{d81}{}
\psfrag{d82}{}
\psfrag{d83}{}
\psfrag{d91}{}
\psfrag{d92}{}
\psfrag{d93}{}
\psfrag{d101}{}
\psfrag{d102}{}
\psfrag{d103}{}
\psfrag{d111}{}
\psfrag{d112}{}
\psfrag{d113}{}
\psfrag{d121}{}
\psfrag{d122}{}
\psfrag{d123}{}
\psfrag{d131}{}
\psfrag{d132}{}
\psfrag{d133}{}
\psfrag{d141}{}
\psfrag{d142}{}
\psfrag{d143}{}
\psfrag{d151}{}
\psfrag{d152}{}
\psfrag{d153}{}
\psfrag{d161}{}
\psfrag{d162}{}
\psfrag{d163}{}
\epsfig{file=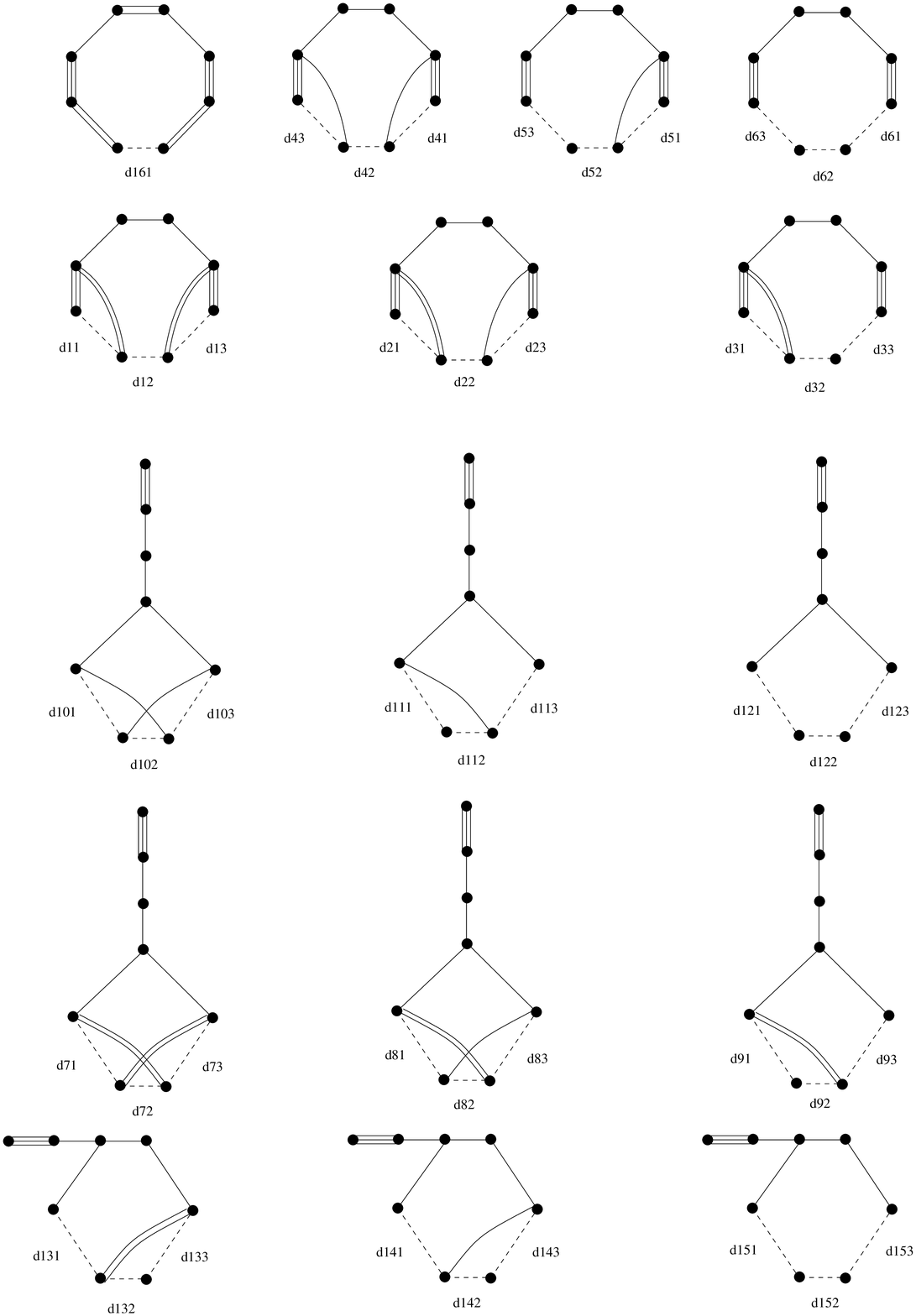,width=0.89179\linewidth}
\end{center}
\end{table}

\pagebreak

\begin{table}[!h]
\caption{$4$-dimensional polytopes with $7$ facets}
\label{t4}
\medskip
\begin{center}
\psfrag{8}{\footnotesize 8}
\psfrag{10}{\footnotesize 10}
\epsfig{file=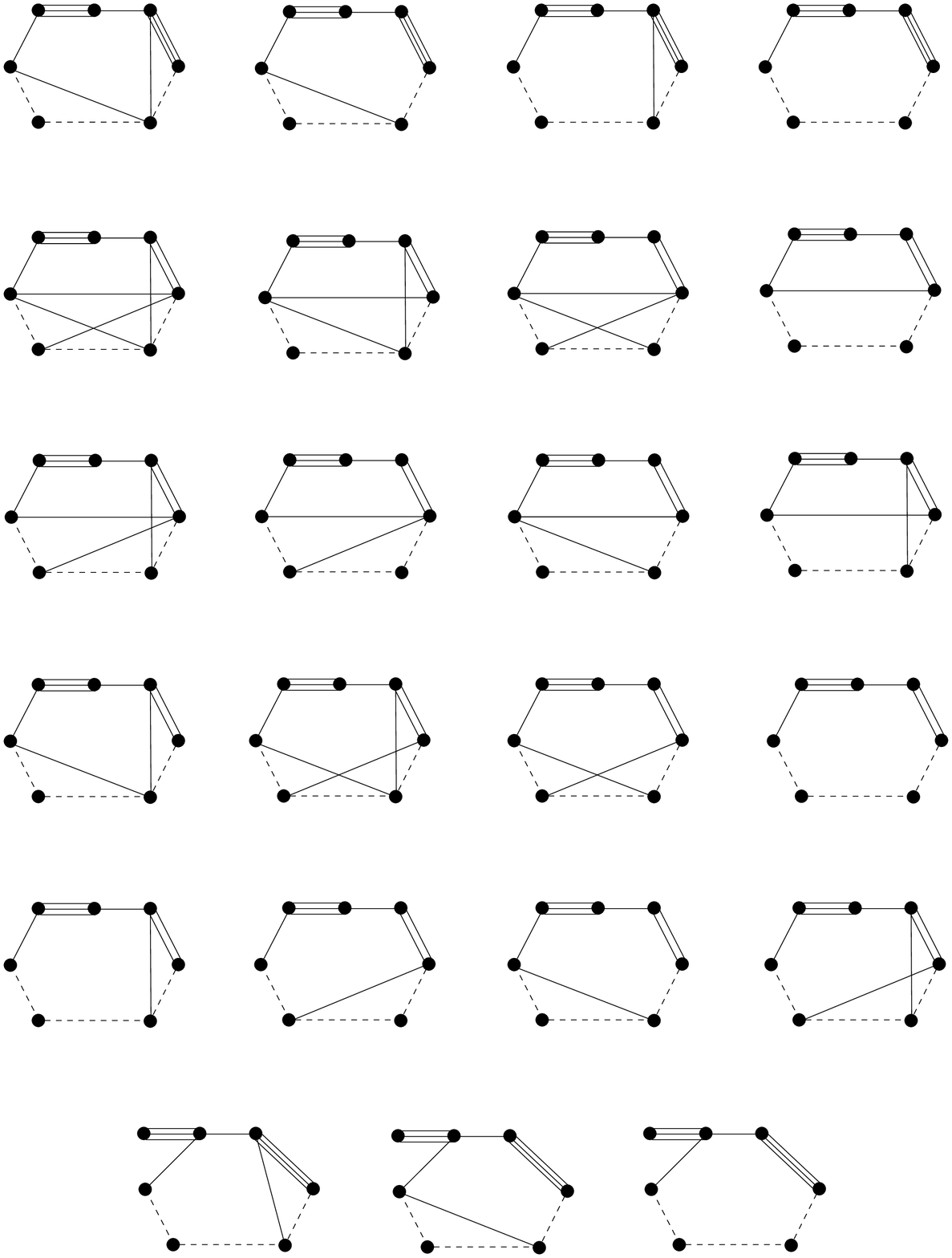,width=0.95\linewidth}
\end{center}
\end{table}

\pagebreak

\begin{center}
{ Table~\ref{t4}: Cont.}

\bigskip
\vspace{20pt}
\psfrag{10}{\footnotesize 10}
\psfrag{8}{\footnotesize 8}
\epsfig{file=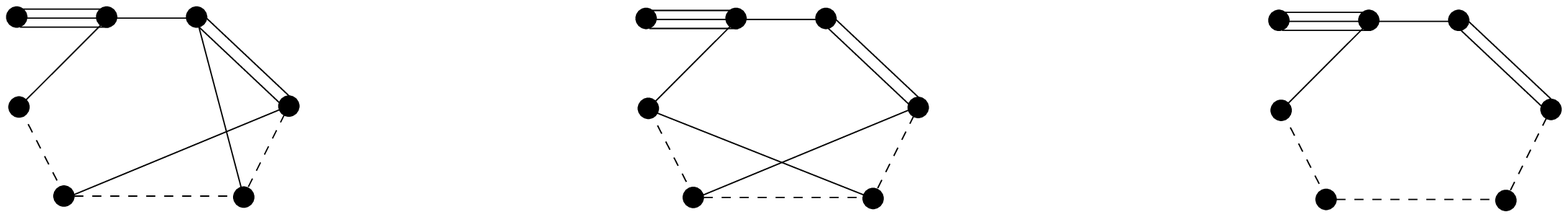,width=0.96\linewidth}\\
\bigskip
\bigskip
\psfrag{10}{\footnotesize 10}
\psfrag{8}{\footnotesize 8}
\epsfig{file=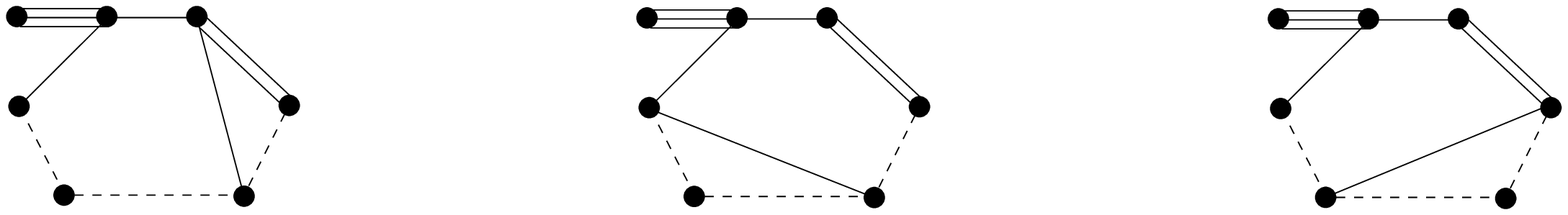,width=0.96\linewidth}\\
\bigskip
\bigskip
\psfrag{10}{\footnotesize 10}
\psfrag{8}{\footnotesize 8}
$\;\;$\epsfig{file=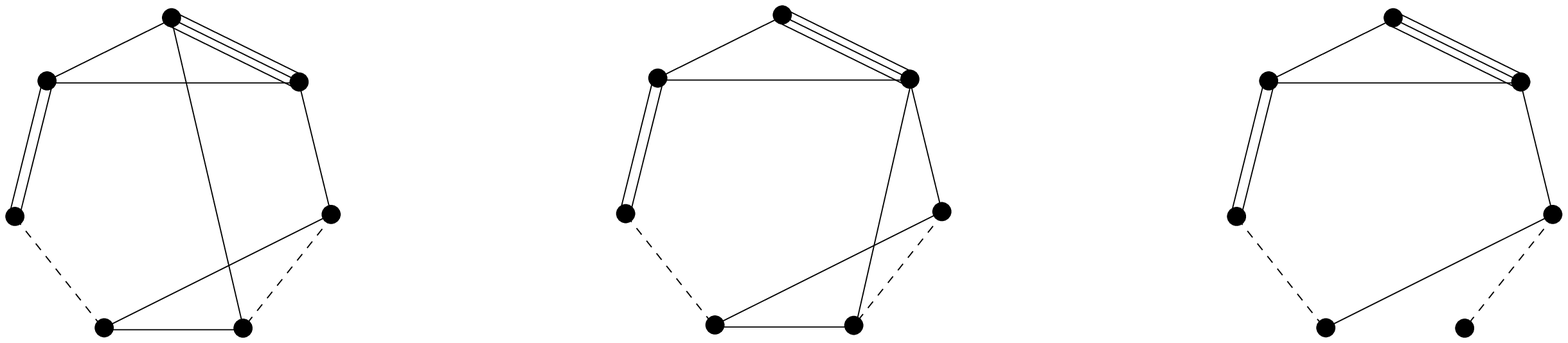,width=0.92\linewidth}\\
\bigskip
\bigskip
\psfrag{10}{\footnotesize 10}
\psfrag{8}{\footnotesize 8}
\epsfig{file=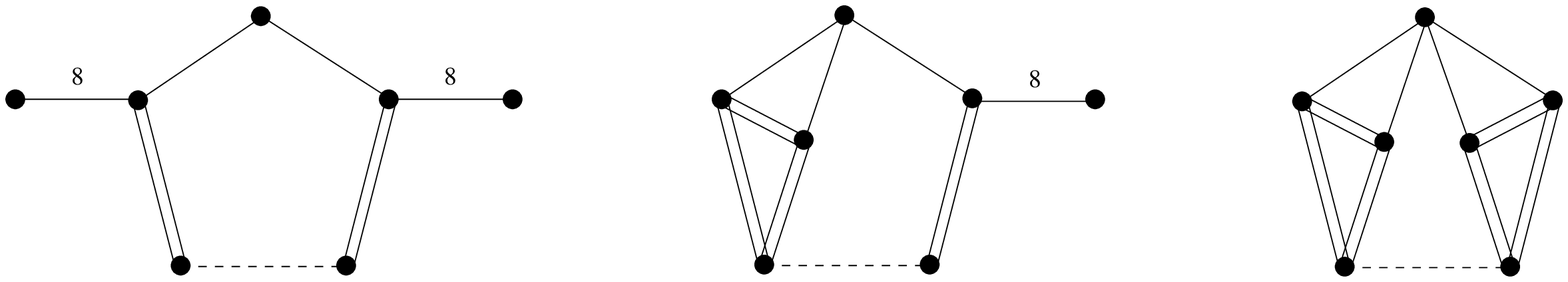,width=0.96\linewidth}\\
\bigskip
\bigskip
\psfrag{10}{\footnotesize 10}
\psfrag{8}{\footnotesize 8}
\epsfig{file=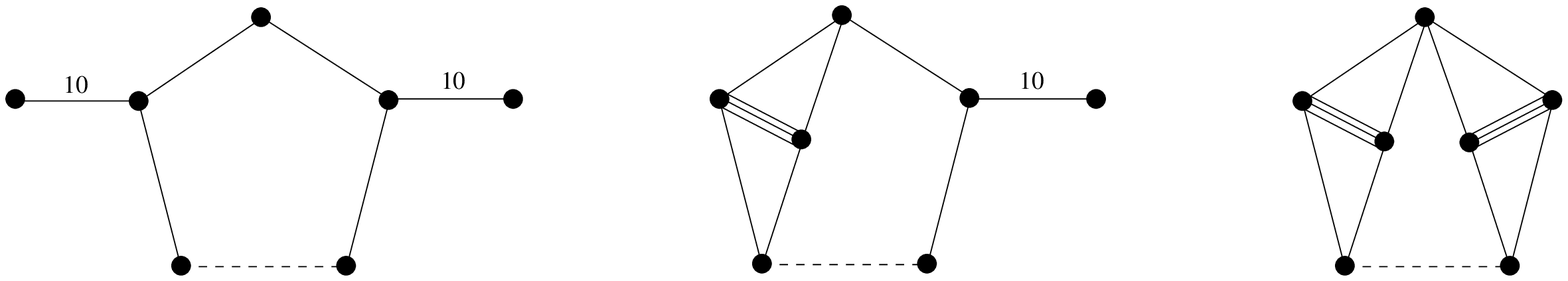,width=0.96\linewidth}\\
\bigskip
\bigskip
\psfrag{10}{\footnotesize 10}
\psfrag{8}{\footnotesize 8}
\epsfig{file=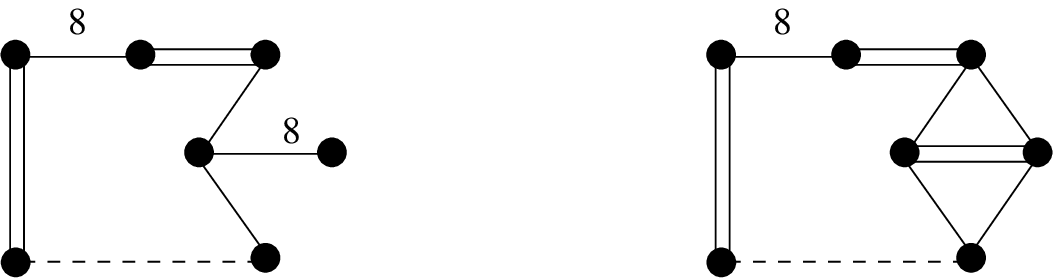,width=0.56\linewidth}
\end{center}


\vfill


\end{document}